\documentclass[11pt]{article}
\usepackage{ amsmath, amssymb,amscd,theorem,chbibref}
\setbibref{\centerline{\large  References}}
{\theorembodyfont{\slshape}\newtheorem{theorem}{Theorem}[section]}
{\theorembodyfont{\slshape}\newtheorem{lemma}[theorem]{Lemma}}
{\theorembodyfont{\slshape}\newtheorem{corollary}[theorem]{Corollary}}
{\theorembodyfont{\slshape}\newtheorem{proposition}[theorem]{Proposition}}
{\theorembodyfont{\slshape}}
{\theorembodyfont{\rmfamily}\newtheorem{example}[theorem]{Example}}
{\theorembodyfont{\rmfamily}\newtheorem{remark}[theorem]{Remark}}
{\theorembodyfont{\rmfamily}\newtheorem{remarks}[theorem]{Remarks}}
\newenvironment{proof}{\noindent{\sc Proof.}}{\vspace{3mm}}
\newenvironment{proofof}{\noindent{\sc Proof}}{\vspace{3mm}}
\numberwithin{equation}{section}

\def \div{\; \big{|} \; }
\def \Aut{\text{\sl Aut}}
\def \Br{\text{\sl Br}}
\def \charr{\text{\sl char}}
\def \cor{\text{\sl cor}}
\def \deg{\text{\sl deg}}
\def \dim{\text{\sl dim}}
\def \expp{\text{\sl exp}}
\def \gcd{\text{\sl gcd}} 
\def \E{ \text{\it \v{E}}}
\def \F{ \text{\it \v{F}}}
\def \Hom{\text{\sl Hom}}
\def \id{\text{\sl id}} 
\def \im{\text{\sl im}} 
\def \ind{\text{\sl ind}} 
\def \inff{\text{\sl inf}} 
\def \kerr{\text{\sl ker}}
\def \L{ \text{\it \v{L}}}
\def \M{ \text{\it \v{M}}}
\def \maxx{\text{\sl max}}
\def \minn{\text{\sl min}}
\def \res{\text{\sl res}}
\def\ov{\overline}
\def\wi{\widetilde}
\def\vg#1{\Delta_{#1}}
\def\indh{\text{\sl{Ind}}_{\ov H\to H}}
\def\zz{\mathbb Z}
\def\si{\sigma}
\def\eig#1#2{{(#1^*/#1^{*p})}^{(#2)}}
\def\GMF{\Gamma_F}
\def\gal{\mathcal G}

\title{\large \bf The First Two Cohomology Groups\\ of
Some Galois Groups}
\date{}
\author{
\normalsize J. Min\'a\v c
\thanks {Supported in part by NSERC grant
R0370A01.}$^{\ \, \dagger}$
\and \normalsize A.
Wadsworth\thanks {The authors would like to thank  
MSRI in Berkeley  for its hospitality while part of the
research for this paper was carried out.  We also thank
the organizers of the exciting Galois theory program
held at MSRI in Fall of 1999. }  }
\begin{document}
\maketitle

\begin{abstract}
We investigate the first two Galois cohomology groups of $p$-extensions
over a base field which does not necessarily contain a primitive $p$th root of unity.
We use twisted coefficients in a systematic way. We describe field extensions which 
are classified by certain residue classes modulo $p^n$th powers of a related field, 
and we obtain transparent proofs
and slight generalizations of some classical results of Albert. The potential application 
to the cyclicity question for division algebras of degree $p$ is outlined. 
\end{abstract}


\setcounter{theorem}{0}
\setcounter{equation}{0}

\setcounter{section}{0}

\section*{\large 
\hfil Introduction\hfil}

\ \indent 
 Let $p$ be a prime number, let $F$ be a field with $\charr(F)\ne p$,
and let $\gal _F
= \gal (F_{sep}/F)$ be the absolute Galois group of~$F$.  It is 
well known that for any $n\in \mathbb N$, the continuous cohomology 
group $H^1(\gal _F, \zz/p^n\zz)$ classifies the cyclic Galois extensions 
of $F$ of degree dividing $p^n$, while $H^1(\gal _F,\mu_{p^n})
\cong F^*/F^{*p^n}$, where $\mu_{p^n}$ denotes the group of $p^n$-th 
roots of unity in the separable closure $F_{sep}$ of $F$; also, 
$H^2(\gal_F, \mu_{p^n}) \cong {}_{p^n}\Br(F)$, the $p^n$-torsion in 
the Brauer group of $F$.  When $\mu_{p^n}\subseteq F$, then 
${\mu_{p^n}\cong \zz/p^n\zz}$ as trivial $\gal _F$-modules, and the
resulting isomorphism $H^1(\gal_F, \zz/p\zz)\cong
H^1(\gal_F, \mu_{p^n})$ is the homological formulation of the classical
Kummer correspondence between cyclic Galois extensions of
$F$ of degree dividing $p^n$ and  cyclic subgroups of $F^*/F^{*p^n}$.
Then also, the Merkurjev-Suslin Theorem describes
${}_{p^n}\Br(F)$; this depends on the isomorphism $H^2(\gal_F,
\mu_{p^n})
\cong H^2(\gal_F, \mu_{p^n}^{\otimes 2})$, which is available because of 
the trivial action of $\gal_F$ on $\mu_{p^n}$.

We consider here questions concerning first and second
cohomology groups with $\mu_{p^n}$ coefficients and concerning 
$F^*/F^{*p^n}$ and ${}_{p^n}\Br(F)$ when $F$ does not contain a primitive 
$p$-th root of unity (so $p\ne 2$), so that the isomorphisms just
mentioned are not available.  In particular, we will give an answer to the 
question:  What field extensions does $F^*/F^{*p^n}$ classify when 
$\mu_{p^n}\not \subseteq F$?  

Since the full absolute  Galois group 
is often too large to work with conveniently, and we are interested in field
extensions of degree a power of $p$, we prefer to work with a
pro-$p$-group instead of $\gal _F$.  For this, let $F(p)$ be the maximal 
$p$-extension of $F$, which is the compositum of all the Galois field
extensions of $F$ of degree a power of $p$, and let $G_F$ be the 
Galois group $\gal(F(p)/F)$;  so, $G_F$ is the maximal pro-$p$ 
homomorphic image of $\gal_F$.   Then, $H^1(G_F, \zz/p\zz) \cong
H^1(\gal_F, \zz/p\zz)$, but \lq\lq $H^1(G_F, \mu_{p^n})$" is undefined
since $\mu_{p^n}\not\subseteq F(p)$.  We will given an interpretation of 
\lq\lq $H^1(G_F, \mu_{p^n})$" in this context, and show that once again it
classifies a certain family of field extensions.  However, these are 
field extensions of $M = F(\mu_{p^n})$, rather than those of $F$.  Indeed, 
our general approach is to relate objects over $F$ to those over 
$L = F(\mu_p)$ and over $M = F(\mu_{p^n})$, since the latter are easier to 
understand because of the presence of enough roots of unity.
Passage from $F$ to $L$ is particularly tractable because 
$p \nmid [L:F]$.

Indications of what happens are provided by Albert's work in \cite{A34}
for the case $n = 1$.  Albert showed that the cyclic degree $p$ field
extensions $S$ of $F$ correspond to certain cyclic degree $p$ extensions
$T$ of $L$.  For, if $T = S\cdot  L$, then $S$ is the unique extension of $F$ 
of degree $p$ within~$T$ (corresponding to the prime-to-$p$ part of the 
abelian Galois group $\gal(T/F)$).  Since $T$ is a $p$-Kummer extension of
$L$,  we have $T = L(\root p \of b)$ for some $b\in L^*$ whose class 
$[b]\in L^*/L^{*p}$ generates the cyclic subgroup associated to $T$ in the 
Kummer correspondence.  The question of classifying cyclic extensions of
$F$ of degree $p$ reduces to determining which cyclic extension of $L$
they generate.  For this, Albert showed that a Kummer $p$-extension 
$T' = L(\root p \of {b'})$ has the form $S'\cdot L$ for some cyclic degree 
$p$ extension $S'$ of $F$ iff $\gal(L/F)$ acts on $[b']$ in $L^*/L^{*p}$ the
same way it acts on $\mu_p$.  A nice way of expressing this is as follows: 
Let $H= \gal (L/F)$.  For each character $\chi\!\!: H \to (\zz/p\zz)^*$
($=$ group of units of the ring $\zz/p\zz$) and each $p$-torsion
$H$-module $A$, there is the $\chi$-eigenmodule of $A$ for the 
$H$ action:  $A^{(\chi)} = \{a\in A\mid h\cdot a = \chi(h) a \text
{ for all }h\in H\}$.  Then, Albert's result can be rephrased:  Cyclic 
$p$-extensions of $F$ correspond to cyclic subgroups of 
$(L^*/L^{*p})^{(\alpha)}$, where $\alpha\!\!:H\to (\zz/p\zz)^*$ is the
cyclotomic character defined by 
$h\cdot \omega = \omega ^{\alpha(h)}$ for each $h\in H$,
$\omega \in \mu_p$. 

We consider here $F^*/F^{*p^n}$ and cyclic Galois extensions of $F$ of
degree $p^n$ for arbitrary $n$.  For this, we work with the field
$M = F(\mu_{p^n})$ instead of $L$.  We show in Cor.~\ref{cor.1.12}
that $F^*/F^{*p^n}\cong (M^*/M^{*p^n})^{\gal (M/F)}$ (which is the 
eigenmodule of $M^*/M^{*p^n}$ for the trivial character of 
$\gal(M/F)$).  We show further in  Th.~\ref{th.1.14} that the cyclic
extensions $K$ of $M$ of degree dividing $p^n$ that correspond to cyclic 
subgroups of $F^*/F^{*p^n}$ are those $K$ which are Galois over $F$
with $\gal(M/F)$ acting on $\gal (K/F)$ by the cyclotomic character
$\alpha\!\! :\gal (M/F)\to (\zz/p^n\zz)^*$ for $\mu_{p^n}$.  In addition, 
when $M = L$, we give in Prop.~\ref{prop.1.7} a small generalization of
Albert's result, by showing that then the cyclic field extensions of $F$ of
degree dividing $p^n$ correspond to the cyclic subgroups of
$(M^*/M^{*p^n})^{(\alpha)}$.  (This correspondence breaks down
whenever $M\ne L$, however---see Remark~\ref{rem.1.8}(a).)

Characters on Galois groups can also be used to define twisted actions 
for their modules:  For any profinite group $G$, any continuous character
$\chi\!\!: G\to (\zz/p\zz)^*$, and any $p^n$-torsion discrete 
$G$-module~$A$, 
define the action of $G$ on $A$ twisted by $\chi$ to be the new action
given by $g*a = \chi(g) (g\cdot a)$ (where $g\cdot a$ denotes the original
$G$-action).  We use such a twisted action in \S 2 to give an interpretation
to \lq\lq $H^i(G_F, \mu_{p^n})$", which as written is not well-defined.  We
define this to mean $H^i(G_F, \wi{\mu_{p^n}})$, where $\wi{\mu_{p^n}}$
denotes $\mu_{p^n}$, but with the action of $\gal = \gal(L(p)/F)$ on it 
twisted by a character $\theta^{-1}$ so that $\gal(L(p)/F(p))$ acts 
trivially on $\wi{\mu_{p^n}}$;  hence, $G_F$ acts on $\wi{\mu_{p^n}}$, 
even though not on $\mu_{p^n}$.  For $H^1$, we give a more specific
interpretation in Th.~\ref{th.2.3}, by showing that 
$H^1(G_F, \wi{\mu_{p^n}})\cong (L^*/L^{*p^n})^{(\theta)}$.  It follows easily 
(see Cor.~\ref{cor.2.5}) that the canonical map $H^1(G_F, \wi{\mu_{p^n}})
\to H^1(G_F, \wi{\mu_{p}})$ is surjective.  This result was needed for the 
paper \cite{MW$_1$}, which was the initial impetus for the work given here.
We also show  that $H^1(G_F, \wi{\mu_{p^n}})$ is isomorphic to an 
eigenmodule of $H^1(G_M, \zz/p^n\zz)$, so that the cyclic subgroups of 
$H^1(G_F, \wi{\mu_{p^n}})$ classify certain cyclic field extensions of $M$ of
degree dividing $p^n$; we give in  
Th.~\ref{th.2.7} a Galois theoretic characterization
of those field extensions of $M$. 

In \S 3 we consider second cohomology groups with $\mu_{p^n} $
coefficients, or, equivalently, the $p^n$-torsion of the Brauer group.
We have the standard isomorphisms ${}_{p^n}\Br(F) \cong
H^2(\gal(L(p)/F), \mu_{p^n})$ and ${}_{p^n}\Br(L) \cong H^2(G_L, 
\mu_{p^n})$.  When $\mu_{p^n}\subseteq L$, the Merkurjev-Suslin 
Theorem gives a very useful presentation of ${}_{p^n}\Br(L)$ by generators
(namely symbol algebras) and relations.  The Merkurjev-Suslin Theorem
does not apply to ${}_{p^n} \Br(F)$ since $\mu_{p^n}\not \subseteq F$, 
but the easy isomorphism ${}_{p^n}\Br(F) \cong 
({}_{p^n}\Br(L))^{\gal(L/F)}$ allows one to analyze ${}_{p^n}\Br(F)$
in terms of the more readily understood ${}_{p^n}\Br(L)$.  This approach
was used by Albert in \cite{A34} in proving his cyclicity criterion for 
algebras of degree $p$, and by Merkurjev in \cite{M83} in proving that 
${}_p\Br(F)$ is generated by algebras of degree $p$.  We give here in 
Th.~\ref{prop.3.6} a generalization of Albert's result, by showing that 
when $\mu_{p^n}\subseteq L$, a division algebra $D$ of degree $p^n$
over $F$ is a cyclic algebra iff there is $d\in D$ with 
$d^{p^n}\in F^*-F^{*p}$.  

But, what we find most tantalizing here is the potential application to the
cyclicity question for division algebras of degree $p$.  If $B$ is a central
division algebra of degree $p$ over$L$ with $[B]\in ({}_p\Br(L))^{\gal(L/F)}$,
then there is a unique central division algebra $A$ over $F$ of degree $p$
with $A\otimes_F L\cong B$.  When $B$ is a cyclic algebra, it is actually a
symbol algebra, i.e., it has a presentation by generators $i,j$ such 
that $i^p = b$, $j^p = c$, $ij = \omega ji$, where $b,c\in L^*$ and 
$\omega \in \mu_p$, $\omega\ne 1$.  We prove in Prop.~\ref{prop.3.4}
that if $A$ is a cyclic algebra, then not only is $B$ cyclic, but it must have 
a presentation as above with $b$ and $c$ mapping to specified 
eigencomponents of $L^*/L^{*p}$ with respect to the $\gal(L/F)$ action.
Thus,  $B$ could very well be cyclic without satisfying the more
stringent conditions which correspond to cyclicity of $A$.  Then $A$ would
be a counterexample to the decades-old question whether all central
simple algebras of degree $p$ must be cyclic algebras (which is 
still unsettled for all $p\ge 5$).  Regrettably, we have not found such a 
counterexample, but we feel that the approach merits further
investigation.

An interesting extreme case of this is for the field $J = 
F(p)\cdot L = F(p)(\mu_p)$.  We have ${}_p\Br(F(p))\cong
 ({}_p\Br(J))^{\gal(J/F(p))}$.  If this group is nonzero, then Merkurjev's 
result says that there is a division algebra of degree $p$ over $F(p)$;
but such an algebra cannot be cyclic, as $F(p)$ has no cyclic field
extensions  of degree $p$.  This observation led us to attempt do find
division algebras of degree $p$ in $({}_p\Br(J))^{\gal(J/F(p))}$, by using
valuation theory.  In \S 4, we describe how valuations on $F$ with residue 
characteristic not $p$ extend to $F(p)$ and to $J$.  This makes it easy to
see that ${}_p\Br(J)$ can have some nontrivial eigencomponents
for the action of $\gal(J/F(p))$ (see Ex.~\ref{laurent}).  But the question 
whether $({}_p\Br(J))^{\gal(J/F(p))} \ne 0$ remains open.

We fix throughout the paper the notation mentioned in this 
Introduction.  The most relevant fields are shown in the diagram below;
the names of the Galois groups given here will also be fixed 
throughout.  We write $\mu_m$ for the group of $m$-th roots of 
unity (in an algebraic closure of the relevant field) and $\mu^*_m$
for the primitive $m$-th roots of unity.  Also, for any profinite group
$G$, let $X(G) = \text{\sl Hom}(G, \mathbb Q/\zz) = H^1(G, \mathbb Q/
\zz)$, the (continuous) character group of $G$.  If $G$ is a group of 
automorphisms of some field, we write $\mathcal F(G)$ for the fixed field 
of $G$.
\parbox{2.5in}{
\hskip1.1truein
\unitlength=1mm
\begin{picture}(22,40)(2,0)
\put(10,-7){\line(4,1){9}}
\put(11,20.25){\line(4,1){8}}
\put(6,-5){\line(0,1){21}} 
\put(22,-2){\line(0,1){5}} 
\put(22,8.5){\line(0,1){11.5}} 
\put(22,26){\line(0,1){8}}
\put(8,25){\line(3,4){8}}
\put(14,-7.75){\makebox(0,0){\scriptsize$H$}}
\put(6,-8){\makebox(0,0){$F$}}
\put(22,-5){\makebox(0,0){$L$}}
\put(22,6){\makebox(0,0){$M$}}
\put(3,5.5){\makebox(0,0){\scriptsize$G_F$}}
\put(6,20){\makebox(0,0){$F(p)$}}
\put(22,23){\makebox(0,0){$J$}}
\put(22,38){\makebox(0,0){$L(p)$}}
\put(10,32){\makebox(0,0){\scriptsize$N$}}
\end{picture}
\bigskip\medskip}
\parbox{3in}{
\addtolength{\jot}{-3pt}
\begin{align*}
\qquad L & = F( \mu_p)\\
M &= F(\mu_{p^n}) \\
J & = F(p) \cdot L = F (p) ( \mu_p)\\
\mathcal G & = \mathcal G (L(p)/F)\\
G_F & = \mathcal G (F(p)/F)\\
G_L & = \mathcal G (L(p)/L)\\
\GMF & = \mathcal G (J/L) \cong G_F\\
N & = \mathcal G (L(p) / F(p))\\
H & = \mathcal G (L/F) \cong \mathcal G(J/F(p))
\end{align*}
 }\hfil

While we were in the last stages of writing this paper we learned of the 
interesting recent preprint by U.~Vishne \cite{V}.  There is some overlap
between this paper and Vishne's (most notably in our Th.~\ref{prop.3.6}),
but it is not great, since Vishne is concerned primarily with the 
situation that $\mu_{p^n}\subseteq F(\mu_p)$.

We would like to thank Bruno Kahn for very helpful discussions at 
an early stage of this work, particularly for pointing out to us
the significance of group actions twisted by characters.

\vskip 0.4truein

\setcounter{section}{1}
\setcounter{theorem}{0}
\setcounter{equation}{0}

\section*{\large 1\ \ Extension fields of $F$ and $L$, and an interpretation
of
$F^* / F^{*p^n}$}

\ \indent
 In this section, we recall some known properties of
$p$-extensions, then give characterizations of the Galois
$p$-extensions of $L = F ( \mu_p)$ which correspond to such
extensions of $F$.  We then look at group actions on
$p^n$-torsion modules, and examine eigenspace
decompositions of such modules.  This is applied to various
structures associated to the abelian $p^n$-extensions of $F$
and $L$.  This leads to  Kummer-like characterization of
$F^* \big / F^{*p^n}$ and $L^* \big / L^{*p^n}$ in terms of certain
abelian
$p^n$-extensions, but the extensions are of
$M = F(\mu_{ p^n})$, not of $F$ and in general not of 
$L$.  
For Prop.~\ref{prop.1.1} through Cor.~\ref{cor.1.3}, $p$ may be any 
prime number.  After Cor.~\ref{cor.1.3} we will assume
further  that $p$ is odd.

We first consider the subfields of $F(p)$.  By definition, $F(p)$ is the 
compositum of all the Galois field extensions of $F$ of degree
a power of $p$ (in some algebraic closure of $F$).  Recall that if 
$K_1$ and $K_2$ are Galois extensions of $F$ with $[K_i:F]=p^{k_i}$,
then $K_1\cdot K_2$ is also Galois over $F$ and $[K_1\cdot K_2:F]
\,\big |\,p^{k_1+k_2}$.  Consequently, 
$F(p)$ can also be described as the 
union of of all the Galois extensions of $F$ of degree a power of $p$.
So, $F(p)$ is the maximal Galois extension of $F$ such that 
$\gal(F(p)/F)$ is a pro-$p$-group.
Also, if $T$ is a field extension of $F$ with $[T:F]<\infty$ and 
$K$ is the normal closure of $T$ over $F$, 
then $T\subseteq F(p)$ iff $T$ is separable over $F$ and 
$[K:F]$ is a power of $p$.
We give next a characterization of such fields $T$ which is
well-known, but we could find no reference for it.
It is an easy
consequence of the property of $p$-groups that every
maximal proper subgroup is a normal subgroup of  index $p$.  
The proof will be omitted.

\begin{proposition}\label{prop.1.1}
Let $S$ be a field of any characteristic, and let $T$ be a field,
$T \supseteq S$, $[T : S ] < \infty$.  Then the following are
equivalent:
\begin{enumerate}
\item[{\rm(i)}]  The normal closure of $T$ over $S$ is Galois
over
$S$ of degree a power of $p$, i.e., $T \subseteq S (p)$.

\item[{\rm(ii)}]  There is a chain of fields $S = S_0 \subseteq
S_1
\subseteq S_2 \subseteq \ldots \subseteq S_k = T$ with each
$S_i$ Galois over $S_{i-1}$, and $[ S_i  : S_{i-1}] = p$.
\end{enumerate}
\end{proposition}

We record two important corollaries of
Prop.~\ref{prop.1.1}, which are also well-known, and
whose easy proofs are omitted.  The second is an immediate
consequence of the first.

\begin{corollary}\label{cor.1.2}  For fields $S$ and $T$ of any
characteristic, if $S \subseteq T \subseteq S(p)$, then $T(p)
= S(p)$.  In particular, $S(p) (p) = S(p)$.
\end{corollary}

\begin{corollary}\label{cor.1.3}  If $T$ is a field with $\charr
(T) \not= p$ and $\mu_p \subseteq T$, then $T(p)^p = T(p)$.
\end{corollary}

From now on, throughout the rest of the paper,  we assume that 
the prime number $p$ is odd.

We now return to the basic setting of this paper, with $F$ a
field, $\charr (F) \not= p$, $L = F ( \mu_p)$, $J = F(p) \cdot L =
F(p) ( \mu_p )$.  
Since $[L:F] \,\big|\, (p-1)$, we have $L\cap F(p) = F$ and
 $\mathcal G (J/L) \cong \mathcal G (F
(p)/F) = G_F$.  Therefore,  there is a canonical inclusion and index
preserving one-to-one correspondence between the
$p$-extensions $S$ of $F$ (i.e., $F \subseteq S \subseteq
F(p)$) and those $p$-extensions $T$ of $L$ with $T \subseteq
J$.  When $S \leftrightarrow T$ we have $T = S\cdot L
S ( \mu_p)$ and
$S = T \cap F(p)$.  
Furthermore, $S$ is Galois over $F$ iff $T$ is Galois over $L$;
when this occurs, $\gal(S/F)\cong \gal(T/L)$.
Of course $L(p)$ may be much larger than
$J$ (as in Ex.~\ref{laurent} below), so we next characterize those Galois
$p$-extensions of $L$ which lie in $J$.  Let $H = \mathcal G
(L/F)$, so $H$ is cyclic and $|H| \div p - 1$.  Let $s = |H| = [L:F]$.    
 
\begin{proposition}\label{prop.1.4}  Let $T$ be a Galois
extension of $L$ with $L \subseteq T \subseteq L(p)$.  Then,
the following are equivalent:
\begin{enumerate}
\item[{\rm(i)}]  $T \subseteq J$.
\item[{\rm(ii)}] $T = S \cdot L$ for a field $S$, $F \subseteq S
\subseteq F(p)$, with $S$ Galois over $F$.
\item[{\rm(iii)}]  $T$ is Galois over $F$ and $\mathcal G (T/F)
\simeq \mathcal G (T/L) \times \mathcal H$ for some group
$\mathcal H$.  (Then necessarily $\mathcal H \cong H$.)
\item[{\rm(iv)}] $T$ is Galois over $F$ and $\mathcal G (T/F)$
has a normal subgroup $\mathcal H$ of order $s$.
\end{enumerate}
If $\mathcal G (T/L)$ is abelian, then (i)--(iv) are equivalent to:
\begin{enumerate}
\item[{\rm(v)}]  $T$ is Galois over $F$ and $\mathcal G (L/F)$ acts 
trivially on $\mathcal G (T/L)$.
\end{enumerate}
\end{proposition}

\begin{proof}
Note that $[T:L]$ could be finite or infinite.
 (i) $ \Leftrightarrow$ (ii) was noted above.  (ii)
$\Rightarrow$ (iii)   Let $\mathcal H = \mathcal G (T/S)$.
Note that 
$S$ and $L$ are each Galois over $F$, and $S \cap L = F$
as $\gcd ([S_0 : F] , [L : F]) = 1$ for each finite degree 
subextension $S_0$ of $F$ in $S$.
Therefore,  $S$ and $L$ are linearly
disjoint over $F$, and $\mathcal G (T/F) = \mathcal G (T/L)
\times \mathcal H$.  Note that $|\mathcal H| = [S\cdot L:S]
=[L:F] = s$.
(iii) $\Rightarrow$ (iv)   For $\mathcal H$
as in (iii), we have $\mathcal H \cong \mathcal G (T/F)
 / \mathcal G (T/L) \cong \mathcal G (L/F) = H$.  Clearly,
$\mathcal H$ is a normal subgroup of $\mathcal G (T/F)$.  
(iv)~$\Rightarrow$~(ii)  
%
%
For $\mathcal  H$ as in (iv), let $S= \mathcal F(\mathcal H)$.  Then
$S$ is Galois over $F$ and $[T:S] = |\mathcal H| = s$.  We have 
$[T:S\cdot L] \, \big|\, [T:S] = s$ and $[T:S\cdot L]$ is a power of 
$p$, as $S\cdot L \subseteq T\subseteq  L(p) = (S\cdot L)(p)$
(see Cor.~\ref{cor.1.2}).  Hence, $S\cdot L= T$.  Since $L$ is 
Galois over $L\cap S$, we have 
$[L:L\cap S] = [L\cdot S:S] = [T:S] = s = [L:F]$;
hence, $L\cap S = F$.  Since $T = S\cdot L$ is a compositum 
of Galois $p$-extensions of $L$ and $\gal(T/L)\cong \gal(S/F)$,
this $S$ is a compositum of such extensions of $F$.  Hence, 
$S\subseteq F(p)$, proving (ii).

Now assume $\mathcal G (T/L)$ is abelian.  Then (iii) implies that
$\mathcal G (T/F)$ is abelian, so (v) holds.  
Conversely, assume (v).  Let $G =
\mathcal G (T/F)$ and let $P = \mathcal G (T/L)$.
Since $G/P$ acts trivially on the abelian group $P$, this 
$P$ must be central in $G$.  Since further $G/P\cong
\gal(L/F)$, which is cyclic, an elementary exercise in group theory
shows that $G$ is abelian.  Let $T_0$ be any finite degree
subextension of $L$ lying in $T$. Then, $[T_0:L] = p^k$, for some
$k\in \mathbb N$, and $T_0$ is abelian Galois over $F$, 
with $|\gal(T_0/F)| = s\kern.08em p^k$.  The primary
decomposition of 
$\gal(T_0/F)$ gives $T_0 = S_0\cdot L_0$, where
$S_0$ is the unique subfield of $T_0$ with $[T_0:S_0] = s$
and $[T_0:L_0] = p^k$. Since $S_0$ is Galois over $F$ with 
$[S_0:F] = p^k$, we have $S\subseteq F(p)$.  Also, 
$[T_0:L] = [T_0:L_0]$, so $L = L_0$.  Thus, $T_0 = 
(T_0\cap F(p))\cdot L$ for every finite degree subextension 
$T_0$ of 
$L$ in $T$.  Hence, $T = (T\cap F(p))\cdot L$, proving (ii).
%
%
  \hfill $\square$
\end{proof}

In proving the following corollary, we will use supernatural numbers.
Recall (cf.~\cite[p.~5]{Se2}) that a supernatural number is a formal 
product 
$ \prod\limits_{i=1}^\infty p_i^{r_i}$ where the $p_i$ are 
distinct prime numbers and each $r_i\in \{0, 1,2, \ldots\}\cup
\{\infty\}$.  Supernatural numbers can be multiplied in the 
obvious way, and likewise the notions of divisibility, gcd's , and 
lcm's of
supernatural numbers have the obvious interpretation.  
For fields $F\subseteq K$ with $K$ algebraic over $F$, we define 
$[K:F]$ to be the supernatural number
$\text{\sl {lcm}} \{\, \dim_FN\mid N\text{ is a field, }  F\subseteq N
\subseteq K,\text { and } \dim_FN<\infty\,\}$.  (Of course, this 
agrees with the usual definition when $ \dim_FN<\infty$.)
The reader can check that for any field $E$ with 
$F\subseteq E\subseteq K$, we have 
\begin{equation}
[K:F] \ =\  [K:E][E:F]\,.
\label{tower}
\end{equation}

\begin{corollary}\label{cor.1.5}
Let $T/F$ be a Galois subextension of $L(p)/F$ such that there 
is a normal subgroup $\mathcal H$ of 
$\mathcal G(T/F)$ with $|\mathcal H| = s = [L:F]$.  Then, 
$L \subseteq T \subseteq J$.  Hence, 
the Galois group $ \mathcal G (L(p) / J)$ is the 
smallest closed normal subgroup $B$ of $\mathcal
G = \mathcal G (L (p)/F)$ such that $\mathcal G / B$ contains a
 normal subgroup of order $s$.    
\end{corollary}

\begin{proof}
We have $[L(p):F] = s\, p^r$ for $0\le r\le\infty$. Since 
$s = |\mathcal H|\div [T:F]$, we have $[L(p):T] = p^t$ for
$t\le \infty$, by (\ref{tower}).  So, $[T(\mu_p):T]$ is a power of
$p$; hence, $\mu_p\subseteq T$, so $L\subseteq T$.  Then, 
Prop.~\ref{prop.1.4} shows that $T\subseteq J$, so of course
$\mathcal G(L(p)/J) \subseteq \mathcal G(L(p)/T)$. 

Now,  $\mathcal G / ( \mathcal G (L(p)/J) \cong \mathcal G (J/F)$,
which contains the normal subgroup $\mathcal G (J/F(p))$ of
order $s$.  The inclusion $T\subseteq J$ just proved 
shows that $\mathcal G (L(p)/J)$ is minimal with this property.
\hfill $\square$
\end{proof}

We want to describe the abelian $p^n$-extensions of $L$ in
$J$ as an eigencomponent of the family of all such extensions
of $L$.  For this, we first need some facts about
$p^n$-torsion $G$-modules.

Fix a  positive integer $n$.  Let $G$ be any profinite group,
and let $A$ be a discrete $G$-module (written additively),
which is $p^n$-torsion as an abelian group.  Let $\chi \!\!  : G
\to (\mathbb Z/p^n \mathbb Z)^*$ be any continuous group
homomorphism.  (The continuity of the character $\chi$ is
equivalent to $\kerr (\chi)$ being an open subgroup of $G$.) 
Note that since $( \mathbb Z /p^n \mathbb Z)^* \cong \Aut (
\mathbb Z / p^n \mathbb Z, +)$ canonically, defining such a 
$\chi$ is equivalent to giving $\mathbb Z /p^n \mathbb Z$ the
structure of a discrete $G$-module.  Because $A$ is
$p^n$-torsion, it is a $\mathbb Z/p^n \mathbb Z$-module.  Let
$A^{(\chi)}$ denote the $\chi$-{\it eigenmodule} of $A$, i.e.
\begin{equation}
A^{(\chi)} 
= \{ a \in A \mid g \cdot a = \chi(g) \cdot a , 
\mbox{ for all } g \in G \} .
\label{eq.1.1}
\end{equation}
Note that $A^{(\chi)} $ is a $G$-submodule of $A$.  If
$A^{(\chi)} = A$, we say that $G$ {\it acts on $A$ via $\chi$}. 
Whatever the original action of $G$ on $A$, we can use $\chi$
to twist the action, obtaining a new $G$-module, denoted
$A_\chi$, such that
$A_\chi = A$ as abelian groups, but if we denote by $\cdot$
the original action of $G$ on $A$ and $*$ the action of $G$ on
$A_\chi$, then 
$$
g * a = \chi (g) (g \cdot a) \qquad
\mbox{for all } g \in G, a \in A .
$$
Another way of saying this is that the map $a \mapsto a
\otimes 1$ is a $G$-module isomorphism
\begin{equation}
A_\chi 
\cong A \otimes_\mathbb Z \mathbb Z / p^n \mathbb Z ,
\label{eq.1.2a}
\end{equation}
where $G$ acts on $\mathbb Z / p^n \mathbb Z$ via $\chi$ (so
$g
\cdot (a \otimes k) = (g \cdot a) \otimes \chi(g) k$).  We will
frequently use the obvious identity
\begin{equation}
(A_\chi)^G
= A^{(\chi^{ -1})} ,
\label{eq.1.2}
\end{equation}
where $\chi^{-1} (g) = \chi (g)^{-1}$.

Now assume $p$ is odd, and suppose $H$ acts on $A$, where $H
= \mathcal G (L/F)$.  Since $A$ is $p^n$-torsion this is
equivalent to $A$ being a module for the group ring $\mathbb
Z /p^n \mathbb Z [H]$.  Now, $H$ is a cyclic group of order~$s$,
where $s \div p - 1$, so 
\begin{equation}
\mathbb Z/p^n \mathbb Z [H]
\cong \mathbb Z / p^n \mathbb Z [x] / 
(x^s - 1 ) .
\label{eq.1.3}
\end{equation}
Since $(\mathbb Z / p^n \mathbb Z)^*$ is cyclic of order
$\varphi (p^n) = (p - 1 ) p^{n-1}$, it contains a cyclic subgroup
of order
$s$.  The elements $\gamma_1 , \dots , \gamma_s$ of this
group are distinct roots of $x^s - 1$, so $x^s - 1 = (x -
\gamma_1) \dots (x - \gamma_s)$ in $\mathbb Z / p^n
\mathbb Z$.  Moreover,  since the map $\mathbb Z / p^n
\mathbb Z \to \mathbb Z / p  \mathbb Z$ has kernel of order
$p^{n-1}$, which is prime to $s$, the subgroup of order $s$
intersects this kernel trivially.  Hence, for $i \not= j$,
$\gamma_i - \gamma_j$ maps to a nonzero element of
$\mathbb Z /p \mathbb Z$; so $\gamma_i - \gamma_j \in
(\mathbb Z/p^n \mathbb Z)^*$.  Therefore, the ideals $(x -
\gamma_i ) \mathbb Z/p^n \mathbb Z [x]$ and $(x -
\gamma_j) \mathbb Z / p^n \mathbb Z [x]$ comaximal
 since
their sum contains the unit $\gamma_i - \gamma_j$.  It
follows by the Chinese Remainder Theorem that
$$
\mathbb Z / p^n \mathbb Z [x] \big/(x^s - 1 )
\;\cong \:\textstyle \bigoplus \limits
^s_{i=1} \mathbb Z / p^n
\mathbb Z [x]  \big/ (x - \gamma_i ) .
$$
Hence, $\mathbb Z / p^n \mathbb Z [H]$ has $s$ mutually
orthogonal primitive idempotents  $e_1, \dots , e_s$, such that
each $e_j \mathbb Z / p^n \mathbb Z [H] \cong \mathbb Z / p^n
\mathbb Z$, and if $h$ is a designated generator of $H$ 
corresponding to the
image of $x$ in \eqref{eq.1.3}, then
$h e_i = \gamma_i e_i$. 
One can check that if $ts \equiv  1\text{ \sl{mod}\ } p^n$, then 
\begin{equation} 
e_i \ = \ t\textstyle\sum\limits_{j=0}^{s-1} \gamma_i^{-j}h^j\,.
\label{eformula}
\end{equation}
If $\chi_i \!\! : H \to
(\mathbb Z / p^n \mathbb Z)^*$ is the character given by $h^j
\mapsto \gamma_i^j$, then in the left multiplicative action on
$e_i \mathbb Z / p^n \mathbb Z [H]$, $H$ acts by $\chi_i$. 
Therefore, for our $H$-module $A$, we have $A =
\bigoplus\limits^s_{i=1} e_i A$, and $H$ acts on $e_i A$ by
$\chi_i$; it follows easily that $e_i A = A^{(\chi_i)}$, and 
\begin{equation}
A = \textstyle \bigoplus\limits^s_{i=1} A^{(\chi_i)} ,
\label{eq.1.4} 
\end{equation}
which is a canonical  decomposition of $A$ into 
a direct sum of eigenmodules.

Now, let $\L$ be the compositum of all the Galois field
extensions of $L$ with cyclic Galois group of order dividing
$p^n$.  Then $\L$ is clearly Galois over $L$, and over $F$, and
$\L$ is the maximal abelian $p^n$-extension of $L$, i.e., the
unique maximal Galois extension of $L$ such that the Galois
group is abelian $p^n$-torsion.  Let $X( \L / L)  = X (\mathcal 
G (\L/L))$ denote the continuous character group of
$\mathcal G ( \L / L)$, 
$$
X( \L / L) = \Hom
(\mathcal G ( \L / L), \mathbb Q / \mathbb Z) \qquad\mbox{ 
(continuous homomorphisms)}.
$$
Note that $X ( \L / L ) \cong \Hom (G_L,  p^{-n} \mathbb
Z / \mathbb Z)$, and that $\mathcal G ( \L / L) \cong G_L
\big/\cap \{ \kerr \, \psi \mid\psi \in \Hom (G_L ,   p^{-n}
\mathbb Z / \mathbb Z) \}$.

Because $\L$ is Galois over $F$ and $\mathcal G ( \L / L)$ is
abelian, there is a well-defined group action of $H =
\mathcal G (L/F)$ on $\mathcal G ( \L/L)$ 
given by,  for $\tau \in \mathcal G
(\L/F)$, $\overline \tau $ the image of $\tau$ in $\mathcal
G(L/F)$, $\sigma \in \mathcal G (\L / L)$, 
$\ov \tau\cdot \sigma = \tau\sigma\tau^{-1}$.
This in turn induces a (left)
action of $H$ on $X( \L / L)$ given by, 
%
%
for $\psi \in X ( \L / L)$,
\begin{equation}
( \overline \tau \cdot \psi ) ( \sigma )
=\psi(\ov\tau^{-1}\cdot \sigma)
= \psi ( \tau^{-1} \sigma \tau ) .
\label{eq.1.5}
\end{equation}
Since $\mathcal G ( \L / L)$ and $X ( \L / L)$ are $p^n$-torsion
$H$-modules, each has an eigenmodule decomposition as
described in \eqref{eq.1.4}: 
\begin{equation}\label{decomp}
\mathcal G ( \L / L)
= \textstyle\prod\limits^s_{i=1} \mathcal G ( \L /
L)^{(\chi_i)}
\qquad
\mbox{and} \qquad
X ( \L / L)
= \textstyle\bigoplus\limits^s_{i=1} X ( \L / L)^{(\chi_i)} .
\end{equation}

These eigendecompositions are related to each other. 
Consider the canonical $\mathbb Z$-bilinear
pairing, 
\begin{equation}\label{pairing}
B \!\!:   \mathcal G (\L / L) \times  X( \L / L) 
\to  p^{-n} \mathbb Z / \mathbb Z \qquad
\mbox{given by } B ( \sigma, \psi) = \psi (\sigma) .
\end{equation}
If $Y$ is any subgroup of $X(\L/L)$, let $Y^\perp = 
\bigcap\limits_{\psi\in Y}\ker(\psi)$ a closed subgroup of 
$\mathcal G (\L/L)$.   Note that $B$ induces an isomorphism 
$Y\cong X(\mathcal G(\L/L)/ Y^{\perp})$.   
Now the pairing in \eqref{pairing} is $H$-equivariant, i.e.,
\begin{equation}
B( h \cdot \sigma, h \cdot \psi)
=h\cdot B ( \sigma, \psi )= B ( \sigma, \psi ), \quad
\mbox{for all } h \in H , 
\sigma \in \mathcal G ( \L / L ), \psi \in X ( \L /L) .
\label{eq.1.6 }
\end{equation}
For any $\chi \in \{ \chi_1, \dots ,\chi_s\}$, let 
\begin{equation}\label {chiL}
{}_\chi L
= \mathcal F(\mathcal N), \text{\ where \ }
\mathcal N = \big (X(\L/L)^{(\chi^{-1})}\big)^\perp\,.
\end{equation}
Then, $\mathcal G(_{\chi}L/L)\cong \mathcal G(\L/L)\big/ \mathcal N$,
hence $X(_{\chi}L/L)\cong X(\L/L)^{(\chi^{-1})}$.  Because 
$H$ acts by $\chi^{-1}$ on $X(_{\chi}L/L)$, it follows from the 
$H$-equivariance of the pairing \eqref{pairing} that $H$ acts by 
$\chi$ on $\mathcal G(_{\chi}L/L)$, i.e., 
$\mathcal G(_{\chi}L/L) = \big(\mathcal G(_{\chi}L/L)\big)^{(\chi)}$.
Moreover, $_{\chi}L$ is the largest subfield of 
$\L$ with this property.

Now, look back at the eigendecompositions of $\mathcal G(\L/L)$ 
and $X(\L/L)$ in \eqref{decomp}.  Because the pairing in \eqref{pairing}
is $H$-equivariant (and $1 - \chi_i(h)\chi_j(h) \in (\zz/p^n\zz)^*$
whenever $\chi_i(h)\chi_j(h) \ne 1$) we have 
 $B(\mathcal G( \L / L)^{(\chi_i)}$,  $X(\L / L)^{(\chi_j)}) =
0$ unless $\chi_j = \chi^{-1}_i $, and $B$ induces a
nondegenerate pairing between $\mathcal G ( \L / L)^{(\chi)} $
and $X ( \L / L)^{(\chi^{-1})}$ for each $\chi \in \{ \chi_1, \dots ,
\chi_s\}$.  Hence, $\big (X(\L/L^{(\chi_i)})\big)^\perp
=\bigoplus \limits_{j\ne i}\mathcal G(\L/L)^{(\chi_j)}$, so that 
%
\begin{equation}
\mathcal G (_\chi L/L) \cong \mathcal G ( \L / L)^{(\chi)}  \qquad
\mbox{ and }  \qquad
\L \cong {}_{\chi_{_1}}\!L \otimes_L \ldots \otimes_L {}_{\chi_{_s}}\!L \, .
\label{eq.1.7}
\end{equation}

\begin{corollary}\label{cor.1.6}
Let $\chi_1$ be the trivial character $H \to ( \mathbb Z / p^n
\mathbb Z)^*$.  Then 
$$
J \cap \L = {}_{\chi_{_1}}\!L \qquad
\mbox{and } \qquad
\mathcal G (( J \cap \L)/L) \cong \mathcal G (\L / L)^H . 
$$
Furthermore, there is a one-to-one correspondence between
the abelian $p^n$-extensions $S$ of $F$ and the abelian
$p^n$-extensions $T$ of $L$ such that $T$ is Galois over $F$
and $H$ acts trivially on $\mathcal G (T/L)$.
\end{corollary}

\begin{proof}
This is immediate from Prop.~\ref{prop.1.4}. \hfill $\square$
\end{proof}

Because $\mu_{p^n} \nsubseteq F$ in general, the abelian
$p^n$-Galois extensions are not classified by subgroups of
$F^* / F^{* p^n}$.  We will show in Th.~\ref{th.1.14} below
that instead  $F^* /
F^{*p^n}$ classifies certain abelian $p^n$-extensions of
$M$, where $M = F ( \mu_{p^n})$. This will require some preliminary 
results. We first make some basic observations about 
$M$ and group actions of $\mathcal G(M/F)$.
Note  that $M \subseteq J$, by Prop.~\ref{prop.1.4}, since
$M$ is Galois over $F$ with $\mathcal G (M/F)$ abelian.

Let $\M$ be the maximal abelian $p^n$-extension of $M$; by
Kummer theory, ${\M = M (\{ \root {p^n} \of m \mid m \in M\})}$.  Then,
 $\gal (M/F)$ acts on the $p^n$-torsion groups 
$\gal (\M/M)$ and $X(\M/M)$  (analogous to the 
action described in \eqref{eq.1.5} above).  So, for each character $\chi
\!\!:\gal(M/F) \to (\zz/p^n\zz)^*$ we have the eigenmodules
$\gal(\M/M)^{(\chi)}$ and $X(\M/M)^{(\chi)}$, though the eigenmodules
together  do not yield a
direct sum decomposition of the whole module, since the 
groupring $\zz/p^n\zz[\gal(M/F)]$ is not semisimple when 
$M\supsetneq L$.  Just as for the  $_\chi L$ defined 
in \eqref{chiL} above, 
we have the field $_\chi M$ which is maximal in 
$\M$ such that $\gal(M/F)$ acts on $\gal(_\chi M/M)$ by 
$\chi$; also $X(_\chi M/M) \cong X(\M/M)^{(\chi^{-1})}$.

%
We can relate these eigenmodules to those of 
%
%
$M^* / M^{* p^n}$. For $m\in M^*$, we write 
$[m]$ for its image in $M^*/M^{*p^n}$.
%
%
From Kummer theory, we have the following isomorphisms,
\begin{equation}
M^*\big  / M^{* p^n} \to \Hom ( \mathcal G ( \M
/M), \mu_{p^n} )
\cong X( \M / M) \otimes_{\mathbb Z} \mu_{p^n} \, ,
\label{eq.1.8}
\end{equation}
where the first map is given by $[m] \mapsto ( \sigma \mapsto
\sigma (a) / a)$, for any $a \in \M^*$ such that $a^{p^n }= m$. 
Now, $\gal(M/F)$ acts on each of the
groups in \eqref{eq.1.8}  (it acts on $\gamma \in
\Hom ( \mathcal G ( \M / M), \mu_{p^n})$ by $(\ov g \cdot
\gamma ) (
\sigma ) = \ov g ( \gamma (g^{-1} \sigma g))$, for any 
$g\in \gal(L(p)/F)$ restricting to $\ov g \in \gal(M/F)$\,);  it is easy to
check that each of the isomorphisms in \eqref{eq.1.8} is
compatible with the $\gal(M/F)$-action.  Thus, if we let
$\alpha$ be the cyclotomic character, which corresponds to the
action of $\gal(M/F)$ on $\mu_{p^n}$, i.e.
\begin{equation}
\alpha \!\! : \gal(M/F) \to ( \mathbb Z/p^n \mathbb Z)^*
\mbox{ is given by } g( \omega ) = \omega^{\alpha (g)}
\mbox{ for all } g \in \mathcal G(M/F)
\mbox { and } \omega \in \mu_{p^n},
\label{eq.1.9}
\end{equation}
then \eqref{eq.1.8} shows that as $\mathcal G$-modules,
\begin{equation}
M^* \big / M^{* p^n} 
\cong X (\M /M)_\alpha  \, .
\label{eq.1.10}
\end{equation}
%
%
It follows that for any character
$\chi\!\!  : H \to ( \mathbb Z /p^n \mathbb Z)^*$,  we have
\begin{equation}
_\chi M= M ( \{ \root p^n\of m 
\mid [m] \in (M^* \big/ M^{*p^n})^{ ( \alpha \chi^{-1})} \}) ,
\label{eq.1.11}  
\end{equation}
since $X ({}_\chi M / M) \cong X ( \M / M)^{(\chi^{-1})} = (X ( \M /
M)_{\alpha} )^{({\alpha} \chi^{-1})}
\cong (M^* \big/ M^{* p^n})^{(\alpha \chi^{-1})}$.
(One can also deduce \eqref{eq.1.11} from the 
$\gal(M/F)$-equivariant Kummer pairing 
$\gal(\M/M)\times M^*/M^{*p^n}\to \mu_{p^n}$.)
Note that if $[m]\in (M^*/M^{*p^n})^{(\alpha\chi^{-1})}$, then 
$M(\root {p^n}\of m)$ is Galois over $F$, since the 
cyclic subgroup $\langle[m]\rangle \subseteq  M^*/M^{*p^n}$
is stable under the action of $\gal(M/F)$.  Conversely, if 
$M(\root {p^n}\of m)$ is Galois over $F$, then $\langle[m]\rangle$
must be $\gal(M/F)$-stable, so $[m]\in (M^*/M^{*p^n})^{(\varphi)}$
for some character $\varphi$.  However, $M(\root {p^n}\of m)$ is
{\it abelian} Galois over $M$ iff $\gal(M/F)$ acts trivially on  
$\gal(M(\root {p^n}\of m)/M)$, iff 
$M(\root {p^n}\of m)\subseteq {}_{\chi_1}\! M$, 
where $\chi_1$ is the trivial character, iff (by 
\eqref{eq.1.11}) $[m]\in (M^*/M^{*p^n})^{(\alpha)}$.

We will relate this to $F^* \big/ F^{*p^n}$ 
below.  But first, let us observe how it yields 
a slight generalization of Albert's 
characterization of the cyclic
 extensions of $F$   of degree $p$ 
(see \cite[Th.~2]{A34} or \cite[p.~211, Th.~15]{A$_1$}).


 \begin{proposition}\label{prop.1.7}
Suppose $p\nmid[F(\mu_{p^n}):F]$.  Let $M = F(\mu_{p^n})$, 
$m\in M^*-M^{*p}$, and $T = M(\root {p^n}\of m)$.  Then, 
$T = S\cdot M$ for some cyclic field extension $S$ of $F$ of degree
$p^n$ iff $T\subseteq J$, 
 iff $g\cdot[m] = [m^{\alpha(g)}]$ in 
$M^*\big/M^{*p^n}$ for each $g \in \mathcal G(M/F)$, 
where $\alpha\!\! :\gal(M/F)\to (\zz/p\zz)^*$
is the cyclotomic character described in \eqref{eq.1.9}
above. 
 \end{proposition}

\begin{proof}
Since $[M:L]$ is a power of $p$, the assumption that 
$p\nmid[F(\mu_{p^n}):F]$ is equivalent to: $M=L$, i.e., 
$\mu_{p^n}\subseteq F(\mu_p)$. For $T = M(\root{p^n} \of m)$, we
have $T = S\cdot M$ for some cyclic Galois extension $S$ of 
$F$ with $[S:F] = [T:M] = p^n$ iff $T\subseteq J\cap \M$  (see 
Prop.~\ref{prop.1.4}).  But, by Cor.~\ref{cor.1.6}, $ J\cap \M=
{}_{\chi_1}M$, where 
$\chi_1\!\!\!: \mathcal G(M/F)\to (\mathbb Z/p^n\mathbb Z)^*$ is the 
trivial character.  As noted in (\ref{eq.1.11}), ${}_{\chi_1}M = 
M(\{\root{p^n}\of c\ | \ [c] \in (M^*\big/M^{*p^n})^{(\alpha)}\,\}$
So, $T = S\cdot M$ iff $T\subseteq {}_{\chi_1}M$ iff (by Kummer theory)
$[m]\in (M^*\big/M^{*p^n})^{(\alpha)}$, as desired.   \hfill $\square$
%
\end{proof}

\begin{remarks}\label{rem.1.8} (a)
Albert's result is the case $n=1$ of Prop.~\ref{prop.1.7}, for which the
assumption that ${p\nmid[F(\mu_p):F]}$ always holds.  Note that the 
hypothesis that $p\nmid[F(\mu_{p^n}):F]$ is needed for 
Prop.~\ref{prop.1.7}.  For, 
%
%
when $M = F(\mu_{p^n})
\supsetneqq F ( \mu_p)$,  we do have $X ( \M / M)^{\mathcal G (M/F)} \cong (
M^* / M^{*p^n})^{(  \alpha )}$ as $\mathcal G (M/F)$-modules (see
\eqref{eq.1.10}).  But the abelian $p^n$-extensions of $M$ coming from such
extensions of $F$ correspond to the image of $X ( \F / F)$ in $X ( \M/M)$, which
is properly smaller than $X ( \M/M)^{\mathcal G (M/F)}$.  Specifically, if we
take $\omega \in \mu_{p^{n^*}} \subseteq M^*$, then clearly $[\omega ] \in
(M^* / M^{*p^n})^{(  \alpha )}$ and $[ \omega ]$ has order $p^n$ in $M^*
/ M^{*p^n}$ (see Lemma~\ref{lem.1.10} below).  Then, for $\psi \in X ( \M / M)$
with $\kerr ( \psi ) = \mathcal G ( \M / M (\root {p^n}\of  \omega))$, we
have $\psi \in X ( \M / M)^{\mathcal G (M/F)}$, but $\psi$ is not in
the image of $X ( \F / F)$, since there is no cyclic
$p^n$-extension $S$ of $F$ with $S \cdot M = M (
\root {p^n}\of  \omega)$. (For, $M(\root p^n \of \omega) \cap F(p)$
is cyclic over $F$, but of degree greater than~$p^n$.)
However, this example is, roughly speaking, the only
obstruction when $M \supsetneqq L$.
For, one can check that $\F \cdot M (\root {p^n}\of  \omega)
={}_{\chi_1}M$, which is  
is the subfield of $\M$ corresponding to $X ( \M
/ M)^{\mathcal G (M/F)}$.  (This follows from the 
observations  that if 
$[M:L] = p^{n-c}$ as in Lemma~\ref{lem.1.10} below, 
then $\big |X(\M/M)^{\mathcal G(M/F)}:\im(X(\F/F))\big| \le
\big |H^2(\mathcal G (M/F), \zz/p^n \zz)\big| = p^{n-c}$ while 
$[\F\cdot M(\root p^n\of \omega) :\F\cdot M] = p^{n-c}$.)

(b)  For a closely related description of the cyclic extensions of $F$ of 
degree $p^n$, see \cite[Th.~2.3]{Sa}.  Saltman's description does not 
require the hypothesis that $\mu_{p^n}\subseteq F(\mu_p)$.
\end{remarks}


We want to relate $F^* \big / F^{*p^n}$ to $M^* \big / M^{*p^n}$, and use 
$F^* \big / F^{*p^n}$
to parametrize certain field extensions of $M$.  For this, we first need the
following two lemmas.

\begin{lemma}\label{lem.1.10}
Let $d = \sup \{ k \mid
\mu_{p^k} \subseteq L \}\in \mathbb N \cup \infty$ and
let $c = \minn (d,n)$.  Then,
$[M : L] = p^{n-c}$.
\end{lemma}

\begin{proof}
We can assume that $n > c$, since otherwise there is nothing
to prove.  Let $\omega$ be a primitive $p^n$-th root of unity,
and let $\nu = \omega^p$ and $\varepsilon =
\omega^{p^{n-c}}\in L$.  Let $M_0 = L( \mu_{p^{n-1}}) = L( \nu)$.  By
induction on $n$ we may assume that $[M_0 : L] = p^{n-c-1}$. 
Then $f(x) = x^{p^{n-c-1}} - \varepsilon$ is the minimal
polynomial of $\nu$ over $L$, since it has the right degree
and $f ( \nu ) = 0$.  Hence, for the norm from $M_0 $ to $L$, we
have $N_{M_0/L} ( \nu) = (-1)^{p^{n-c-1}} ( - \varepsilon ) =
\varepsilon$ (as $p$ is odd).  

Since $M = M_0 ( \omega ) = M_0 ( \nu ^{1/p})$ and $\mu_p
\subseteq M_0$, $M$ is a $p$-Kummer extension of $M_0$.  So
$[M : M_0] = p $ or $=1$.  If $[M:M_0] = 1$, then $\omega \in
M_0$, so $\varepsilon = (N_{M_0 /L} ( \omega ))^p$, so
$N_{M_0 / L} ( \omega)$ is a primitive $p^{c+1}$ root of
unity in $L$, contradicting the definition of $c$ (as $n > c$). 
Thus, we must have $[M : M_0 ] = p$, so $[M:L] = [M :
M_0][M_0 : L] = p^{n-c}$ as desired.  \hfill $\square$ 
\end{proof}

\begin{lemma}\label{lem.1.11}
$H^1 ( \mathcal G(M/L), \mu_{p^n} )
= H^2 ( \mathcal G(M/L), \mu_{p^n} ) = 1$.
\end{lemma}

\begin{proof}
Because $\mathcal G (M/L)$ is cyclic (as $p$ is odd), $H^2 (
\mathcal G (M/L), \mu_{p^n}) \cong ( \mu_{p^n} )^{
\mathcal G (M/L)} \big / N_{M/L} (\mu_{p^n} )$.  Now $(
\mu_{p^n})^{\mathcal G (M/L)} = \mu_{p^n} \cap L =
\mu_{p^c}$.  Let $\omega \in
\mu_{p^{n^*}}$, and let $\varepsilon = \omega^{p^{n-c}} \in
\mu_{p^{c^*}} \subseteq L$.  Since $M = L( \omega )$ and $[M :
L] = p^{n-c}$ by Lemma~\ref{lem.1.10}, $\omega$ has minimal
polynomial $x^{p^{n-c}} - \varepsilon$ over $L$, so $N_{M/L}
( \omega ) = (-1)^{p^{n-c}} ( - \varepsilon ) = \varepsilon$. 
Thus, $N_{M/L} (\mu_{p^n} ) = \langle N_{M/L} ( \omega )
\rangle = \langle \varepsilon \rangle = \mu_{p^c} = (
\mu_{p^n}) ^{\mathcal G (M/L)}$.  Hence $H^2 ( \mathcal G
(M/L),
\mu_{p^n }) = 1$, and, as $\mathcal G (M/L)$ is cyclic and
$\mu_{p^n}$ is finite, the Herbrand quotient \cite [p.~134, Prop.~8]{Se}
shows that
$\big |H^1 (
\mathcal G ( M/L), \mu_{p^n} )\big | =\big |H^2 ( \mathcal G (M/L),
\mu_{p^n})\big |$.  \hfill $\square$
\end{proof}

\begin{corollary}\label{cor.1.12} \ 
\begin{enumerate}
\item[\rm{(i)}]  $L^* \big /L^{* p^n} \cong (M^* \big /
M^{*p^n})^{\mathcal G (M/L)}$.

\item[\rm{(ii)}]  $F^* \big / F^{*p^n} \cong (L^*
\big /L^{*p^n})^{\mathcal G (L/F)} \cong (M^* \big /M^{*p^n})^{\mathcal
G (M/F)}$.
\end{enumerate}
\end{corollary}

\begin{proof}
(i) From the $5$-term exact sequence of low
degree terms associated with the Hochschild-Serre spectral
sequence (cf.~\cite[p.~307, Th.~11.5]{R}, the following is exact: 
$$
H^1 ( \mathcal G (M/L) , \mu^{G_M}_{p^n})
\to H^1 (G_L, \mu_{p^n}) 
 \to H^1 (G_M , \mu_{p^n})^{G(M/L)}
\to H^2 (\mathcal G (M/L), \mu_{p^n}^{G_M}) .
$$
Since $\mu^{G_M}_{p^n}= \mu_{p^n}$, Lemma~\ref{lem.1.11}
applies, yielding 
\begin{equation}
H^1 (G_L, \mu_{p^n} )
\cong H^1 (G_M, \mu_{p^n})^{\mathcal G ( M/L)} .
\label{eq.1.14}
\end{equation}
The long exact cohomology sequence arising from the short exact
sequence of $G_L$-modules
\begin{equation}
\begin{CD}
1 @>>> \mu_{p^n} @>>> L(p)^*
@> ( \ )^{p^n} >> L(p)^* @>>> 1
\end{CD}
\label{eq.1.14a}
\end{equation}
(the $p^n$-power map $L(p) \to L(p)$ is onto by
Cor.~\ref{cor.1.3}), yields $H^1 (G_L, \mu_{p^n}) \cong L^*
\big /L^{*p^n}$; likewise $H^1 (G_M, \mu_{p^n} ) \cong M^* \big /
M^{*p^n}$.  Thus, the isomorphism in \eqref{eq.1.14}
translates to $L^*\big  / L^{*p^n} \cong{ (M^* \big /M^{*p^n})^{\mathcal
G(M/L)}}$, as desired.

(ii)   The first isomorphism in (ii) can be proved
in the same way as (i), using that \break ${H^i ( \mathcal G (L/F), (
\mu_{p^n})^{G_L} ) = 1}$ for $i = 1, 2$; this is clear, as $\gcd ( |
\mathcal G (L/F) | , | (\mu_{p^n})^{G_L}|) = 1$.  But the desired
isomorphism can also be obtained easily nonhomologically: 
Injectivity of the map $F^* \big / F^{*p^n} \to (L^*
\big /L^{*p^n})^{\mathcal G(L/F)}$ follows by an evident norm
argument; surjectivity of this map follows using Hilbert's
Th.~90.  The second isomorphism in (ii) follows from the
 isomorphism in part (i).  \hfill
$\square$
\end{proof}

\begin{remarks}\label{rem.1.13} (a)  Lemma~\ref{lem.1.11} and
Cor.~\ref{cor.1.12} are special to $n$ where $M = L (
\mu_{p^n})$.  One can compute analogously to
Lemma~\ref{lem.1.11} that for $k \le n$ and any $i \ge 1$ 
$$
\big | H^i ( \mathcal G (M/L), \mu_{p^k})\big  | = |\mu_{p^a} / \mu_{p^b}|
= p^{a-b} ,
$$
where $a = \minn (k,c)$ and $b = \maxx (k + c - n , 0)$.  So, for
$0 < k < n$, $\big |H^i ( \mathcal G (M/L), \mu_{p^k})\big | > 0$ and the
map $L^* \big/ L^{*p^k} \to (M^* \big / M^{*p^k})^{\mathcal G (M/L)}$ is
neither $1$-$1$ nor onto.

(b)  On the other hand, for every $k\ge n$, the map 
$L^*\big /L^{*p^k} \to M^*\big /M^{*p^k}$ is injective.  More 
generally, for any field $K \supseteq L$, if the map 
$L^*\big /L^{*p^n} \to K^*\big /K^{*p^n}$ is injective, 
then the map $L^*\big /L^{*p^k} \to K^*\big /K^{*p^k}$ is injective 
for every $k\ge n$.   This is a special case of the following 
group theoretic observation:  If $A\subseteq B$ are abelian groups 
such that the map $A/p^nA \to B/p^n B $ is injective and if 
$p^{k-2}({}_{p^{k-1}}B)\subseteq A$ for all $k>n$, then the 
map  $A/p^kA \to B/p^k B$ is injective for all $k>n$.  
(Here ${}_{p^{k-1}}B$ denotes the $p^{k-1}$-torsion subgroup of $B$.)
This is easily proved by induction on $k$.  
One can also check that  $L^*\big/L^{*p^k} \cong 
(M^*\big /M^{*p^k})^{\mathcal G(M/L)}$ for every $k\ge n$.
\end{remarks}

We can now give an answer to the question:  What do $F^*
\big/ F^{*p^n}$ and $L^* \big/ L^{*p^n}$ classify when $\mu_p
\nsubseteq F$?  The answer, a kind of generalized Kummer
theory, is given not in terms of field extensions of $F$ or $L$,
but those of $M = F ( \mu_{p^n})$.  In the following 
Theorem~\ref{th.1.14} we are no longer assuming that $\mu_p
\not\subseteq F$.  Hence, the $F$ appearing there could be 
either the $F$ or the $L$ of the preceding results.

\begin{theorem}\label{th.1.14}
Let $p$ be an odd prime number.  Let $F$ be any field with
$\charr (F) \not= p$, and let $M = F ( \mu_{p^n})$.   Let
$\alpha \!\!: \mathcal G (M/F)
\to ( \mathbb Z/p^n
\mathbb Z)^*$ be the 
cyclotomic character, as in
\eqref{eq.1.9}.  Then,
there is a one-to-one correspondence
between the subgroups $U$ of $F^* / F^{*p^n}$ and the abelian
Galois extensions $K$  of $M$ such that $\expp
( \mathcal G (K/M))  \div  p^n$, $K$ is Galois over
$F$, and
$\mathcal G (M/F)$ acts on $\mathcal G (K/M)$ via $\alpha$.
\end{theorem}

\begin{proof}
One lemma needed for this proof is deferred until the end of this section
for ease of exposition.
Let $\mathcal G =
\mathcal G (F ( \mu_p)(p) / F)$, and let $\alpha$ denote also 
the composition  $\mathcal G \to  \mathcal G (M/F)
\overset{ \alpha}{\longrightarrow} (
\mathbb Z/p^n
\mathbb Z)^*$, which is the cyclotomic
character for the action of $\mathcal
G$ on
$\mu_{p^n}$.

 We have $\mathcal G$-module isomorphisms
\begin{equation}
M^*\big / M^{*p^n} \to \Hom ( \mathcal G ( \M / M), \mu_{p^n})
\to X ( \M /M) \otimes \mu_{p^n}
\to X( \M /M)_\alpha ,
\label{eq.1.15} 
\end{equation}
where the first map is the one in Kummer theory:  $[m] \to
(\sigma \mapsto \sigma (a) / a)$ for any $a \in \M$ with
$a^{p^n} = m$.  The last map depends on a choice of
generator  $\omega$ of $\mu_{p^n}$, and is given by $\psi \otimes
\omega^j \mapsto j \psi$.  By composing these isomorphisms
with the identity map $\iota \!\! : X ( \M /M)_\alpha \to X( \M / M)$
(not a $\mathcal G$-homomorphism), we have a bijective
function
$f \!\! : M^* \big / M^{*p^n} \to X ( \M /M)$.  By Cor.~\ref{cor.1.12} (i) or
(ii), the isomorphisms of \eqref{eq.1.15} map $F^* \big/ F^{*p^n}$
onto $(X ( \M / \M)_\alpha)^{\mathcal G (M/F)}$ which $\iota$
sends onto $X(
\M / M)^{(\alpha^{-1})}$.  Thus, $f$ yields a one-to-one
correspondence between subgroups of $F^* \big / F^{*p^n}$ and
subgroups of $X(\M / M)^{(\alpha^{-1})}$.  

Now, the subgroups $Y$ of $X(\M/M)$ are in one-to-one
correspondence with the closed subgroups $C$ of $\mathcal G (
\M / M)$, where $Y = \{ \psi \in X ( \M /M) \mid \kerr (\psi )
\supseteq C\}$ and $C = \bigcap\limits_{\psi \in Y} \kerr ( \psi
)$.  In the Galois correspondence, every closed subgroup $C$ of
$\mathcal G ( \M / M)$ corresponds to an intermediate field $K$,
$M
\subseteq K \subseteq \M$, where $K = \mathcal F (C)$, the fixed
field of $C$, and $C = \mathcal G ( \M / K) \cong G_K / G_{\M}$. 
Suppose $Y \leftrightarrow C \leftrightarrow K$ (so $Y = X (K/M)
\subseteq X ( \M / M)\,$).  We have seen that $Y$ corresponds to a
subgroup of $F^*\big / F^{* p^n}$ iff $Y \subseteq X ( \M /M)^{(
\alpha^{-1})}$; by Lemma~\ref{lem.1.9}  below (with $G = \mathcal G$, $P
= G_{\M}$, $Q = G_M$, $R = G_K$), this occurs iff $G_K$ is normal
in $\mathcal G$ (i.e., $K$ is Galois over $F$) and $\mathcal G$ acts
on $\mathcal G (K/M)$ via $\alpha$.  The last condition is
equivalent to:  $\mathcal G (M/F)$ acts on $\mathcal G(K/M)$ via
$ \alpha$.
  \hfill
$\square$
\end{proof}

\begin{remarks}\label{rem.1.15}
(a) In the correspondence of Th.~\ref{th.1.14}, a subgroup $U$
of $F^* \big / F^{*p^n}$ corresponds to the field $K = M (\{
\root {p^n} \of b \mid [b] \in U\})$.  Conversely, for a given field $K$,
the corresponding subgroup of $F^* \big / F^{*p^n}$ is $\{ [b] \mid
b \in K^{p^n} \cap F^*\}$.  As usual in Kummer Theory, the
correspondence is a lattice isomorphism, so it preserves
inclusions, intersection and joins;  when
$ U
\leftrightarrow K$ we have $|U| = [K : M]$, and $U$ and
$\mathcal  G(K/M)$ are (Pontrjagin) dual to each other.  In
particular, $|U| < \infty$ iff $|\mathcal G (K/M)| < \infty$, and
when this occurs, $U \cong \mathcal G (K/M)$ (noncanonically).
 
(b)  One would prefer a correspondence between cyclic
subgroups of $F^* \big / F^{*p^n}$ and simple radical extensions
of $F$ of degree dividing $p^n$.  But this does not work when
$\mu_{p^n} \nsubseteq F$.  For example take any $b \in F^*$
with $[b]$ of order $p^n$ in $F^*\big  / F^{* p^n}$.  Let $K = F (a)$ for
some choice of $a$ with $a^{p^n} = b$.  Then (as $b \notin
F^{*p}$ with $p$ odd, so $x^{p^n} - b$ is irreducible in $F[x]$),
$[K:F] = p^n$, and also $[K \cdot M : M] = p^n$, by
Cor.~\ref{cor.1.12}.  Hence, $K$ and $M$ are linearly disjoint
over $F$.  So, when $M \not= F$, i.e., $\mu_{p^n} \nsubseteq F$,
then $\mu_{p^n} \nsubseteq K$, so $K$ does not contain all
the $p^n$-th roots of $b$.  Thus, the field $K$ depends on the
choice of $a$ among the $p^n$-th roots of $b$, but $K \cdot M$
does not depend on such a choice.  For another example, for
$L = F ( \mu_p)$, assume $\mu_{p^n} \nsubseteq L$ and let $\varepsilon \in
\mu^*_{p^c}
\subseteq L$ with $c$ as large as possible, as in Lemma~\ref{lem.1.10}. 
Then $[L ( \root {p^n} \of\varepsilon) : L] = p^n$, by
Lemma~\ref{lem.1.10}, but $[M ( \root {p^n} \of\varepsilon) : M] =
p^c$, so $[ \varepsilon ]$ has order $p^c$ in $L^* \big / L^{*p^n}$ by
Cor~\ref{cor.1.12}.

(c)  Theorem~\ref{th.1.14} indicates the divergence between
$p^n$-th power classes and cyclic $p^n$ field extensions of
$F$ when $\mu_{p^n} \nsubseteq F$:  We have $F^* \big / F^{*p^n}
\cong (M^* \big / M^{*p^n})^{\mathcal G(M/F)}$, while cyclic field
extensions of $F$ of degree dividing $p^n$ correspond to
cyclic subgroups of $X(\F / F)$, which map (not quite
isomorphically) to $X ( \M /M)^{\mathcal G(M/F)} \cong (H^1
(G_M , \mu_{p^n} )_{ \alpha^{-1}})^{\mathcal G(M/F)}
\cong (M^*\big  /M^{*p^n} )^{(  \alpha) }$.  So the $p^n$-th
power classes of
$F$ and cyclic extensions of $F$ correspond to different
eigenspaces for the $\mathcal G (M/F)$ action on $M^* \big /
M^{*p^n}$.

(d)  For Th.~\ref{th.1.14} we required that $p$ be odd.  For $p
= 2$ the theorem is valid, with the same proof, for any field
$F$ $( \charr (F) \not= 2)$ such that $\mu_4 \subseteq F$, and $M =
F(\mu_{2^n})$ for any $n$.  But if $\mu_4 \nsubseteq F$, the
theorem fails already for $n = 2$, since then $-4 \notin
F^{*4}$ but $-4 \in F(\mu_4)^{*4}$ (as $(1+\sqrt{-1})^4 = -4$),
so the map $F^* \big / F^{*4}
\to F(\mu_4)^* \big / F(\mu_4)^{*4}$ is not injective.
\end{remarks}

The following lemma will complete the proof of 
Theorem \ref{th.1.14}.
 For the lemma, let $G$ be a profinite group and let
$P$ and $Q$ be closed normal subgroups of
$G$ with $P \subseteq Q$.  Suppose $Q/P$ is an abelian torsion
group of exponent dividing some $e \in \mathbb N$.  Let $\chi \!\!
: G \to (\mathbb Z / e \mathbb Z)^*$ be any continuous
group homomorphism. 

\begin{lemma}\label{lem.1.9}
Let $R$ be a closed subgroup of $G$ with
$P \subseteq R \subseteq Q$.  Then, $X (Q/R)
\subseteq X (Q/P)^{(\chi)}$ iff $R$ is normal in $G$ and $G$ acts on
$Q/R$ via $\chi^{-1}$.
\end{lemma}

\begin{proof} The action of $G$ on $Q$ by conjugation induces the
action of $G$ on $Q/P$; then $G$ acts on $X(Q/P)$ by $(g
\cdot
\psi) (qP) = \psi (g^{-1} \cdot (qP)) = \psi (g^{-1} q g)P)$, for all
$g \in G$,  $\psi \in X (Q/P)$, $q \in Q$.   For every closed
subgroup
$R$ of $G$ with $P \subseteq R \subseteq Q$, the character group
$X (Q/R)$ embeds in $X (Q/P)$ using the surjection $Q/P \to Q/R$;
we use this embedding to view $X (Q/R) \subseteq X (Q/P)$.

Now, suppose $X (Q/R) \subseteq X (Q/P)^{(\chi)}$.  This means
that
\begin{equation}
\psi^\prime (g^{-1} \cdot (qP))
= \psi^\prime (q^{\chi (g)} P) \qquad
\mbox{for all}  \qquad
q\in Q , \ g \in G,\  \psi \in X (Q/R) ,
\label{eq.1.12} 
\end{equation}
where $\psi^\prime$ denotes the image of $\psi$ in $X(Q/P)$. 
Observe that 
$\bigcap\limits_{\psi \in X(Q/R)}
\kerr (\psi^\prime) = R/P$. It follows from (\ref{eq.1.12})
that for each $r\in R$ and $g\in G$ we have 
$(g^{-1}rg)P = r^{\chi(g)}P\in R/P$.  Hence, $R$ is a normal 
subgroup of $G$.
So, $G$ acts on $Q/R$, and \eqref{eq.1.12}
translates to
\begin{equation}
\psi (g^{-1} \cdot (qR))
= \psi (q^{\chi (g)} R) \qquad
\mbox{for all}  \qquad
q \in Q , \  g \in G, \ \psi \in X (Q/R) .
\label{eq.1.13} 
\end{equation}
Because $\bigcap\limits_{\psi \in X (Q/R)} \kerr ( \psi )$ is trivial,
it follows that $g \cdot (qR) = (qR)^{\chi (g^{-1})}$ for all $q \in Q$,
$g
\in G$, i.e.,
$G$ acts on
$Q/R$ via $\chi^{-1}$, as desired.

Conversely, if $R$ is normal in $G$ and $G$ acts on $Q/R$ via $\chi^{-1}$, then
$g \cdot (qR) = (qR)^{\chi^{-1} (g)}$ for all $g \in G$, $q \in Q$; so
\eqref{eq.1.13} holds, and this shows that $G$ acts on $X (Q/R)$ via $\chi$.
\hfill $\square$
\end{proof}

\vskip 0.4truein

\setcounter{section}{2}
\setcounter{theorem}{0}
\setcounter{equation}{0}

\section*{\large 2\ \ ``$H^1 (G_F, \mu_{p^n})$''}

\ \indent
In the analysis of Demushkin groups as Galois groups in
\cite{MW$_1$} and \cite{MW$_2$}, the authors needed to show that
for any field $F$ ($\charr (F) \not= p$) there is an action of $G_F$ on
$\mathbb Z/p^n \mathbb Z$ so that the map $H^1 (G_F , \mathbb
Z/p^n \mathbb Z) \to H^1 (G_F , \mathbb Z / p \mathbb Z)$ is
surjective.  This does not hold in general for the trivial action on
$\mathbb Z / p^n \mathbb Z$ (see Remark~\ref{rem.2.6} below)
but the authors pointed out in \cite[proof of Th.~2.2]{MW$_2$} that
this does hold if we replace $\mathbb Z / p^n \mathbb Z$ by
$\mu_{p^n}$.  Of course $\mathbb Z / p^n \mathbb Z \cong
\mu_{p^n}$ as abstract groups, but they  have different
$G_F$-actions, and the action of $G_F$ on $\mu_{p^n}$ is defined
iff $\mu_{p^n} \subseteq F(p)$, iff $\mu_p \subseteq F$.  The
present paper was originally motivated by the need to handle the
case where $\mu_p \nsubseteq F$ (so $p$ is odd).    We will do this
in Cor.~\ref{cor.2.5} below, by realizing ``$H^1 (G_F , \mu_{p^n})$'' as
an eigenmodule of $L^* / L^{* p^n}$, where $L = F ( \mu_p)$ (see
Th.~\ref{th.2.3}).  We also show in Th.~\ref{th.2.7} that the
subgroups of ``$H^1 (G_F , \mu_{p^n})$'' classify certain abelian
$p^n$-extensions of $F( \mu_{p^n})$ and such extensions of 
$F(\mu_{p^n})\cap F(p)$.

We assume throughout this section that $p$ is an odd prime
number, and that $F$ is a field with $\charr (F) \not= p$ and that
$\mu_p \nsubseteq F$.  Then let $L = F ( \mu_p)$, and let $J$, $G_F$,
$\mathcal G$, and $H$ be as defined in the Introduction.  Let $\GMF
= \mathcal G (J/L)$.  We need to give meaning to ``$H^1 (G_F ,
\mu_{p^n})$'', since $\mu_{p^n}$ is not a $G_F$-module.  One
approach would be to take this to mean $H^1 (\GMF , \mu_{p^n})$
since $\GMF \cong G_F$ and $\GMF$ does act on $\mu_{p^n}$.  But we
will take a somewhat different approach by replacing $\mu_{p^n}$
by a twisted version $\widetilde{\mu_{p^n}}$ on which $G_F$ does
act; then ``$H^1(G_F, \mu_{p^n})$'' will be replaced by $H^1 (
G_F , \widetilde{\mu_{p^n}})$, which is in fact isomorphic to $H^1
(\GMF , \mu_{p^n})$---see Th.~\ref{th.2.3}(iii).  We will need to
keep track of how twisting a module by a character affects the
action of the Galois group on the cohomology groups of a normal
subgroup, and this is described by our first lemma.   

\begin{lemma}\label{lem.2.1}
Let $G$ be a profinite group, and let $K$ be a closed normal
subgroup of $G$.  Let $A$ be a discrete $p^n$-torsion
$G$-module, and let $\chi \!\! : G \to (\mathbb Z / p^n \mathbb
Z)^*$ be a continuous character such that $K \subseteq \kerr
(\chi)$.  Then, for the $\chi$-twisted $G$-module $A_\chi$ as in
\eqref{eq.1.2a} above,
$$
H^i (K, A_\chi )
\cong H^i (K, A)_\chi
$$
as $G$-(i.e., $G/K $-) modules.
\end{lemma}

\begin{proof}
Let $\cdot$ denote the action of $G$ on $A$ and 
on $Z^i(K, A)$. 
Let $*$ denote the action of $G$ on $A_\chi$ and 
on $Z^i(K, A)_\chi$. 
So for the identity map $j \!\! : A_\chi \to
A$ given by $a \mapsto a$, we have $j (g * a) = \chi (g) (g \cdot
j (a))$ for all $g \in G$, $a \in A$.  Since $j$ is a $K$-module
isomorphism, it induces a group isomorphism on continuous
cocycles, $j^* \!\! : Z^i  (K, A_\chi) \to Z^i (K, A)$, given by $j^*
(f) = j \circ f$.  For $k \in K$, $g \in G$, write $k^g$ for $g^{-1}
kg$.  Then, for all $g \in G$, $f \in Z^i (K, A_\chi)$, $k_1, \dots ,
k_i \in K$, we have
\begin{align*}
[ j^* (g \cdot f)] (k_1, \dots, k_i)
& = j(g * f(k^g_1, \dots, k^g_i )) \\
&= \chi(g) (g \cdot j (f (k^g_1, \dots , k^g_i )))\\
&= \chi (g) [ ( g \cdot j^* (f))(k_1 , \dots , k_i )] .
\end{align*}
So, $j^* (g \cdot f) = \chi (g) g \cdot j^* (f) = g * j^* (f)$,
showing that $j^*$, viewed as a (bijective) mapping $Z^i (K,
A_\chi )$ to $Z^i (K, A)_\chi$ is a $G$-module homomorphism. 
Thus, $j^*$ induces a $G$-module isomorphism $H^i (K, A_\chi )
\to H^i (K, A)_\chi$.  \hfill $\square$
\end{proof}

We now return to our specific setting where $\mathcal G =
\mathcal G (L(p)/F)$, $N = \mathcal G (L(p)/F(p))$, etc.

\begin{lemma}\label{lem.2.2}
Let $A$ be any $p$-primary torsion abelian group which is a
discrete $\mathcal G$-module on which $N$ acts trivially (i.e.,
$A$ is a $G_F$-module).  Then, there are $\mathcal G$-module
maps
\begin{align*}
H^1 (G_F, A) 
& \cong H^1 (\GMF, A)
\cong H^1 (G_L, A)^H , \qquad \mbox{and}\\
H^2 (G_F, A) 
& \cong H^2 (\GMF, A)
\hookrightarrow H^2 (G_L, A).
\end{align*}
\end{lemma}

\begin{proof}
Consider the commutative diagram of $\mathcal G$-module
maps
\begin{equation}
\begin{CD}
H^i (G_F, A) @> \inff>> & H^i ( \mathcal G, A)\\
@VVV & @VV\res V\\
H^i (\GMF, A) @> \inff>> & H^i (  G_L, A)
\end{CD}\label{eq.2.1}
\end{equation}
Here the left vertical map is the isomorphism induced by the
isomorphism $G_F \to \GMF$, which is compatible with the
actions of these groups on $A$, and with the action of $\mathcal
G$ on these groups.  The right vertical restriction map in
\eqref{eq.2.1} has image lying in $H^1(G_L, A)^{\mathcal G} =
H^1 (G_L, A)^H$ (as $H \cong  \mathcal G /G_L)$, and $G_L$ acts
trivially on $H^i(G_L , -)\,$).  Indeed, this map is injective with
image all of $H^i (G_L, A)^H$ as $A$ is $p$-primary and $G_L$ is
the $p$-Sylow subgroup of $\mathcal G$.  (For, the map
$\cor \circ \res \!\! :H^i ( \mathcal G, A) \to H^i (
\mathcal G,A)$ is multiplication by $| \mathcal G / G_L | = |H|$
on the $p$-primary group $H^i ( \mathcal G , A)$, so it is an
isomorphism.  This shows $\res$ is injective.  On the other hand,
$\res \circ \cor \!\! : H^i (G_L, A) \to H^i (G_L, A)$
is multiplication by $|H|$ on $H^i (G_L, A)^H$, so it maps this
$p$-primary subgroup to itself.)

Now, recall the five term exact sequence of low degree terms
on cohomology associated with the short exact sequence $1
\to N \to \mathcal G \to G_F \to 1$:
\begin{equation}
\minCDarrowwidth 8 pt  
\begin{CD}
0 @>>> H^1 (G_F, A^N) @> \inff >> H^1 (\mathcal G, A)
@> \res>> H^1 (N, A)^{G_F} @>>>
H^2 (G_F, A^N) @> \inff >> H^2 ( \mathcal G, A)
\end{CD}\label{eq.2.2}
\end{equation}
Since $N$ acts trivially on $A$, we have $H^1 (N, A) = \Hom (N, A)$
(continuous homomorphisms); but $\Hom (N,A) = 0$ as $A$ is
$p$-primary torsion and $N$ has no cyclic factor groups of
order $p$ (as $F(p)$ is $p$-closed).  Thus \eqref{eq.2.2} shows that 
the inflation map $H^1 (G_F, A) \to H^1 ( \mathcal G, A)$ is an
isomorphism, hence \eqref{eq.2.1} shows $H^1 (\GMF, A) \to
H^1 (G_L, A)^H$ is an isomorphism.  Similarly, \eqref{eq.2.2}
shows that $H^2 (G_F , A^N) \to H^2 (\mathcal G , A)$ is
injective, so
\eqref{eq.2.1} yields the injectivity of $H^2 ( \GMF, A) \to H^2
(G_L, A)$.  \hfill $\square$
\end{proof}

We now apply these lemmas to the $\mathcal G$-module
$\mu_{p^n}$.  We need the following  
characters:  First let $\alpha \!\! : \mathcal G
\to ( \mathbb Z/p^n \mathbb Z)^*$ be the cyclotomic
character for the action of $\mathcal G$ on
$\mu_{p^n}$; that is, for all $\omega \in
\mu_{p^n}$ and all $g \in \mathcal G$
\begin{equation}
g ( \omega ) = \omega^{\alpha (g)} ,
\label{eq.2.3}
\end{equation}  
analogous to  \eqref{eq.1.9} above.  Then, let $\theta \!\! :
\mathcal G \to ( \mathbb Z /p^n \mathbb Z)^*$
be the unique character such that 
\begin{equation}
\theta |_N =\alpha \qquad \mbox{and} \qquad
\theta|_{G_L} = 1.                            
\label{eq.2.4}
\end{equation} 
(Such a $\theta$ is unique because $\mathcal G= N\cdot \mathcal G_L$
and it exists since $\alpha |_{N\cap G_L} = 
\alpha|_{\mathcal G(L(p)/J)} = 1$.)
Indeed, observe that $\theta = \alpha^{p^{n-1}}$ has the
properties specified in \eqref{eq.2.4}.
For, since $( \mathbb Z/p^n
\mathbb Z )^*
\cong ( \mathbb Z/p^{n-1}  \mathbb Z ) \times ( \mathbb
Z/(p-1) \mathbb Z)$, our $\alpha$ must map the pro-$p$ group
$G_L$ into $\mathbb Z/p^{n-1} \mathbb Z$, so $\alpha^{p^{n-1}}|_{G_L}
=1$;  also, as $| N / G_J| = |\mathcal H |$
and $|\mathcal H| \div  ( p - 1 )$,
this $\alpha$ maps $N$ into $\mathbb
Z/(p-1) \mathbb Z$, so $\alpha^{p^{n-1}}|_N = \alpha|_N$.
 Note, in fact that $\theta$ coincides with
the prime-to-$p$ part of $\alpha$, i.e., the composition of
$\alpha$ with the projection $( \mathbb Z / p^n \mathbb Z)^* \to
\mathbb Z/ (p-1) \mathbb Z$.  Likewise, $\alpha \theta^{-1}$ is
the $p$-primary part of $\alpha$.  To get a $G_F$-module
structure from $\mu_{p^n}$, we need to twist $\mu_{p^n}$ by a
character to get a trivial $N$-action.  But we use $\theta^{-1}$
rather than $\alpha^{-1}$ in order to have a character trivial on
$G_L$, so that we can invoke Lemma~\ref{lem.2.1}.

\begin{theorem}\label{th.2.3}  Let $\widetilde{\mu_{p^n}} = (
\mu_{p^n} )_{\theta^{-1}}$ for the character $\theta$ defined
in \eqref{eq.2.4} above.  Then  $\widetilde{\mu_{p^n}}$ is a
$G_F$-module, and 
\begin{enumerate}
\item[{\rm(i)}] $H^1 (G_F, \widetilde{\mu_{p^n}}) \cong H^1
(G_L,
\mu_{p^n})^{(\theta)} \cong (L^* / L^{*p^n})^{(\theta)} $;

\item[{\rm(ii)}] $H^2 (G_F, \widetilde{\mu_{p^n}})
\hookrightarrow H^2 (G_L,
\mu_{p^n})^{(\theta)}$;

\item[{\rm(iii)}]  For all $i \ge 0$, $H^i (G_F,
\widetilde{\mu_{p^n}} )\cong H^i (\GMF,
\mu_{p^n})$ as groups, but $\mathcal G$ acts trivially on $H^i
(G_F , \widetilde{\mu_{p^n}})$, while $\mathcal G$ acts via
$\theta$ on $H^1 (\GMF, \mu_{p^n})$.
\end{enumerate}
\end{theorem}

\begin{proof}
Since $N$ acts on $\mu_{p^n}$ via $\alpha$ and $\theta|_N =
\alpha$, $N$ acts trivially on the $\theta^{-1}$-twist
$\widetilde {\mu_{p^n}}$ of $\mu_{p^n}$.  Hence, there is a
well-defined action of $G_F \cong \mathcal G/N$ on $\widetilde
{\mu_{p^n}}$.  Because $G_F \cong \GMF$, with the isomorphism
compatible with the action of $\mathcal G$ on these groups, and
with the action of these groups on $\widetilde {\mu_{p^n}}$, we
have the $\mathcal G$-module isomorphism $H^i (G_F ,
\widetilde {\mu_{p^n}} ) \cong H^i (\GMF , \widetilde
{\mu_{p^n}})$, for all $i
\ge 0$.  To prove (i), we have the isomorphisms
\begin{align}
H^1 (G_F, \widetilde {\mu_{p^n}} )
&\cong H^1 (\GMF, \widetilde {\mu_{p^n}} )
\cong H^1 (G_L , \widetilde {\mu_{p^n}} )^H \notag\\
&\cong (H^1 (G_L , \mu_{p^n} )_{\theta^{-1}} )^H
= H^1 (G_L , \mu_{p^n} )^{( \theta )}
\cong (L^* / L^{* p^n} )^{(\theta )} ,
\label{eq.2.5}
\end{align}
where the first isomorphism was just noted, the second is by
Lemma~\ref{lem.2.2}, the third by Lemma~\ref{lem.2.1} (since
$\theta|_{G_L}$ is trivial).  For (ii), likewise, we have $H^2 (G_F ,
\widetilde {\mu_{p^n}} ) \hookrightarrow H^2 (G_L, \widetilde
{\mu_{p^n}})^H \cong H^2 (G_L , \mu_{p^n})^{(\theta )}$.

(iii)  We have noted already that $H^i (G_F, \widetilde
{\mu_{p^n}}) \cong H^i (\GMF, \widetilde {\mu_{p^n}} )$ as
$\mathcal G$-modules.  Since $\theta|_{G_L} = 1$, we have
$\widetilde {\mu_{p^n}} \cong \mu_{p^n}$ as $\GMF$-modules, so
$H^i (\GMF , \widetilde {\mu_{p^n}}) \cong H^i (\GMF, \mu_{p^n})$
as abelian groups.  It remains to check the $\mathcal
G$-actions.  For this, note that as $\mathcal H = \mathcal G
(J/F(p))$ acts trivially on $\GMF$ and trivially on $\widetilde
{\mu_{p^n}}$, it must act trivially on $H^i (\GMF , \widetilde
{\mu_{p^n}})$.  Since $\GMF$ also acts trivially on $H^i (\GMF,
\widetilde {\mu_{p^n}} )$, $\mathcal G(J/F) ( \cong \GMF \times
\mathcal H)$ acts trivially on $H^i(\GMF, \widetilde {\mu_{p^n}}
)$.  Because the $\mathcal G$-action on $H^i  (\GMF,
\widetilde{\mu_{p^n}})$ is via $\mathcal 
G (J/F)$, this $\mathcal G$-action is also trivial; so
$\mathcal G$ acts trivially also on
$H^i (G_F ,
\widetilde {\mu_{p^n}} )$.  On the other hand, the fact that
$\mathcal G$ acts trivially on $H^i (\GMF , \widetilde {\mu_{p^n}}
) \cong H^i (\GMF , \mu_{p^n})_{\theta^{-1}}$ translates to
$\mathcal G$ acts via $\theta$ on $H^i (\GMF, \mu_{p^n})$.  \hfill
$\square$
\end{proof}

\begin{remark}\label{rem.2.4}
Since $\widetilde{\mu_{p^n}} ^{\otimes j} = (
\mu_{p^n}{}^{\otimes j} )_{\theta^{-j}}$, the analogue to
Th.~\ref{th.2.3} holds for every $j
\in \mathbb Z$ (with the same proof), with $\widetilde{
\mu_{p^n}}^{\otimes j}$ replacing $\widetilde{\mu_{p^n}}$ and
$\theta^j$ replacing $\theta$, except for the second isomorphism in
\ref{th.2.3}(i). 
\end{remark}

We can now prove the result needed for \cite{MW$_1$}.

\begin{corollary}\label{cor.2.5}
For any $n \in N$, the map $H^1 (G_F , \widetilde{\mu_{p^n}}) \to
H^1 ( G_F , \widetilde{\mu_p})$ induced by the canonical
epimorphism $\widetilde{\mu_{p^n}} \to \widetilde {\mu_p}$ is
surjective.
\end{corollary}

\begin{proof}
Let $\alpha' \!\! : \mathcal G \to ( \mathbb Z / p \mathbb Z )^*$ be
the cyclotomic character given by $\sigma ( \omega )
=\omega^{\alpha^\prime ( \sigma )}$ for all $\sigma \in \mathcal
G$ and $\omega \in \mu_p$, i.e., $\alpha^\prime$ is the $n = 1$
version of $\alpha$; let $\theta^\prime = \alpha^\prime$.  Note
that $\alpha^\prime$ is the composition of $\alpha \!\!: \mathcal G
\to (\mathbb Z / p^n \mathbb Z)^*$ with the canonical
epimorphism $(\mathbb Z / p^n \mathbb Z)^* \to ( \mathbb Z /p
\mathbb Z)^*$; also $\theta^\prime$ is obtainable likewise from
$\theta$.  Therefore, for any $p$-torsion $\mathcal G$-module $A$,
twisting the $\mathcal G$-action on $A$ by $\theta^\prime$ is the
same as twisting by $\theta$.  This applies to $A = \mu_p$ and to
$A = L^*/L^{* p}$; in particular, the two possible interpretations  of
$\widetilde{\mu_p}$  coincide.

Now, consider the diagram
\begin{equation}
\begin{CD}
H^1 (G_F , \widetilde{\mu_{p^n}}) @> \cong >>( L^* /
L^{*p^n})^{(\theta)}\\
@VVV @VVV\\
H^1 (G_F , \widetilde{\mu_ p }) @> \cong >>( L^* /
L^{*p})^{(\theta^\prime)}
\end{CD}\label{eq.2.6}
\end{equation}
 The horizontal maps here are the isomorphisms of
Th.~\ref{th.2.3}(i) for $n$ and for $1$.  By checking the maps that
led to the isomorphisms (see\eqref{eq.2.5}), we can see that
diagram \eqref{eq.2.6} is commutative.  Note also that as $\theta
|_{G_L} = 1$, $\theta$ may be  viewed as a character on $H$.  Hence
$(L^* / L^{* p^n})^{(\theta)}$ is one of the $H$-eigenmodules in
the direct decomposition of $L^* / L^{*p^n}$ as in \eqref{eq.1.3}
above, and $(L^* / L^{*p})^{( \theta^\prime)}$ is the
corresponding
$H$-eigenmodule of $L^* /L^{*p}$.  Since the epimorphism $L^* /
L^{*p^n} \to L^* / L^{*p}$ is an $H$-module map, the right vertical
map in \eqref{eq.2.6} is onto.  Hence the left map in \eqref{eq.2.6} is
also onto, as desired.  \hfill $\square$
\end{proof}

\begin{remark}\label{rem.2.6}
The surjectivity proved in Cor.~\ref{cor.2.5} definitely depends on
the choice of the group action on the cyclic group of order $p^n$. 
Note, by contrast, that for the case of trivial $G_F$-action the
canonical map
\begin{equation}
H^1 (G_F, \mathbb Z/p^n \mathbb Z)
\longrightarrow H^1 (G_F, \mathbb Z/p  \mathbb Z)
\label{eq.2.7}
\end{equation}
is not in general onto.  For, surjectivity of the map in \eqref{eq.2.7}
is equivalent to:  every cyclic Galois extension of $F$ of degree $p$
embeds in a cyclic Galois extension of degree $p^n$.  This is not
always true, as the following example illustrates:  Let $S$ be the
rational function  field 
$S = \mathbb  Q (x_1, \dots, x_p)$ and let $\sigma$ be
the $\mathbb Q$-automorphism of $S$ given by 
cyclically permuting the
indeterminates.  Let $F$ be the fixed field of $\sigma$; so $S$ is
Galois over $F$ of degree $p$.  Let $L = F ( \omega)$ where
$\omega \in \mu^*_p$, and let $T = S( \omega ) = 
\mathbb Q( \omega ) (
x_1, \dots , x_p)$, which is a cyclic Galois extension of $L$ of
degree $p$.  If $T$ lies in a cyclic Galois extension $K$ of $L$ of
degree $p^2$, then a theorem of Albert \cite[p.~207, Th.~11]{A$_1$} says
that there is $\alpha \in T$ such that $N_{T/L} ( \alpha ) = \omega$. 
However, using the unique factorization in $\mathbb Q( \omega )
[x_1,
\dots, , x_n]$, one sees that then there is a constant $\beta \in
\mathbb Q ( \omega)$ such that $N_{T/L} (\beta) = \omega$.  But
$N_{T/L} (\beta) = \beta^p$ as $\beta \in L$, and $\omega$ is
clearly not a
$p$-th power in $T$.  So, there can be no such field $K$, and hence
there is no cyclic extension of $F$ of degree $p^2$ containing $S$.
\end{remark}

The subgroups of $H^1 (G_F, \widetilde {\mu_{p^n}} )$ classify
certain field extensions.  To describe this, let $M = F(
\mu_{p^n})$, and let $E = M \cap F( p)$, which is a cyclic Galois
field extension of $F$ with $E ( \mu_p) = M$.  Slightly abusing notation,
let ${\alpha \theta^{-1}}$ denote also the character for
$\mathcal G(E/F)$ and for $\mathcal G (M/F)$ induced by $\alpha
\theta^{-1}$ on $\mathcal G$;likewise,  let $\alpha$ denote also the
character for $\mathcal G (M/F)$ induced by $\alpha$ on $\mathcal G$.

\begin{theorem}\label{th.2.7}
The subgroups of $H^1 (G_F , \widetilde {\mu_{p^n}})
$ are in one-to-one correspondence with the abelian $p^n$
field extensions $S$ of $E$ such that $S$ is Galois over $F$ and
$\mathcal G(E/F)$ acts on $\mathcal G (S/E)$ via $
{\alpha\theta^{-1}}$.  They are also in one-to-one
correspondence with the abelian $p^n$ field extensions $T$ of
$M$ with $T$ Galois over $F$, $\, T \subseteq J$, and $\mathcal G
(M/L)$ acts on $\mathcal G (T/M)$ via $
\alpha$.
\end{theorem}

\begin{proof}
We have
\begin{align*}
H^1 (G_F , \widetilde {\mu_{p^n}} )
& \cong H^1 (G_L, \mu_{p^n})^{(\theta )}
\cong H^1 (G_M, \mu_{p^n})^{(\theta )} \cong H^1 (G_M , (
\mathbb Z/p^n \mathbb Z)_\alpha )^{(
\theta)}\\
& \cong ( H^1 (G_M , ( \mathbb Z/p^n \mathbb Z)_\alpha )^{(
\theta)}
=  H^1 (G_M ,  \mathbb Z/p^n \mathbb Z)^{(
\theta \alpha^{-1})}
= X( \M / M)^{(\theta \alpha^{-1})} ,
\end{align*}
where the first isomorphism is by Th.~\ref{th.2.3}(i), the
second is the restriction of the $\mathcal G$-isomorphism $H^1
(G_L , \mu_{p^n}) \cong H^1 (G_K, \mu_{p^n})^{\mathcal G
(M/L)}$ (see \eqref{eq.1.14}), the third is from the $\mathcal
G$-isomorphism $\mu_{p^n} \cong ( \mathbb Z / p^n \mathbb
Z)_\alpha$, the fourth by the $\mathcal G$-isomorphism in
Lemma~\ref{lem.2.1}.  By Lemma~\ref{lem.1.9}, the subgroups of
$X ( \M / M)^{(\theta \alpha ^{-1})}$ are in one-to-one
correspondence with the fields $T$ such that $M \subseteq T
\subseteq \M$ (i.e., $T$ is an abelian $p^n$-extension of $M$),
 $T$ is Galois over $F$, and $\mathcal G$ 
(hence also $\mathcal G(M/F)$)
acts on $\mathcal G
(T/M)$ via $\alpha \theta^{-1}$.
When
this occurs, 
$\mathcal G (M/E)$, which  lies in $\kerr ( 
{\alpha \theta^{-1}})$,  acts trivially on $\mathcal G (T/M)$, so
by Prop.~\ref{prop.1.4} (with $E$ replacing $F$), $T \subseteq
J$.  Also, since $G_L \subseteq \kerr ( \theta)$, the 
character  ${\alpha \theta^{-1}}$ agrees with $\alpha$ 
on $\mathcal G(M/L)$.  So, the conditions on $T$ are the ones
stated in the theorem.  Conversely, if $T$ satisfies these
conditions, then since $\mathcal G (M/L)$ acts on $\mathcal
G(T/M)$ via $\alpha$ (which agrees with 
 ${\alpha \theta^{-1}}$ on $\mathcal G(M/L)$)
and $\mathcal G(T/E)$ acts via ${\alpha \theta^{-1}}$
(i.e., trivially) on $\mathcal G (T/M)$, we have $\mathcal G
(M/F)$ acts on $\mathcal G(T/M)$ via ${\alpha
\theta^{-1}}$, so $T$ is involved in the one-to-one
correspondence.  These fields $T$, since $T \subseteq J$,
correspond by Prop.~\ref{prop.1.4} to fields $S \subseteq F(p)$
(where $S = T \cap F(p)$ and $T = S ( \mu_p)\,$) such that $S$ is
Galois over $F$, $\mathcal G (S/E) \cong \mathcal G (T/M)$, and
$\mathcal G (E/F)$ acts on $\mathcal G(S/E)$ the same way as
$\mathcal G (M/L)$ acts on $\mathcal G (T/M )$, i.e., via
${\alpha \theta^{-1}}$.  \hfill $\square$
\end{proof}

\begin{remark}\label{rem.2.8}
An alternate way of obtaining the first one-to-one
correspondence of Th.~\ref{th.2.7} is to observe that
\begin{align*}
H^1(G_F , \widetilde{\mu_{p^n}}) 
&
\cong H^1
(G_E , \widetilde{\mu_{p^n}})^{\mathcal G(E/F)} 
\cong (H^1
(G_E , \mathbb Z / p^n \mathbb Z)_{\alpha
\theta^{-1}})^{\mathcal G(E/F)}
 \\
&
= H^1 (G_E , \mathbb Z / p^n \mathbb Z)^{(\theta \alpha^{-1})}
= X ( \E / E )^{( \theta \alpha^{-1})}\, ,
\end{align*}
where the first
isomorphism follows as in the proof of Cor.~\ref{cor.1.12}(i),
since $H^i ( \mathcal G (E/F) , \widetilde{\mu_{p^n}}) \cong H^i (
\mathcal G (M/L), \widetilde{\mu_{p^n}} ) \cong H^i ( \mathcal G
(M/L) , \mu_{p^n} ) = 0$ for $i = 1, 2 $, by
Lemma~\ref{lem.1.10}.  The subgroups of $X( \E / E)^{(\theta
\alpha^{-1})}$ correspond to the specified abelian $p^n$-extensions of $E$ by Lemma~\ref{lem.1.9}.
\end{remark}

\vskip0.4truein

\setcounter{equation}{0}
\setcounter{theorem}{0}
\setcounter{section}{3}

\section*{\large 3\ \ Cyclic algebras of degree $p$}

\ \indent 
One of the oldest unsettled questions in the theory of central
simple algebras is whether every division algebra of prime degree
$p$ must be a cyclic algebra.   This is known to be true for $p = 2$
and~$3$, but unsettled for every $p \ge 5$.  We now describe a
possible approach to finding a counterexample which arises from
the analysis of the relations between structures over $F$ and
those over $L = F ( \mu_p)$.  We have not succeeded in using this
approach to obtain a counterexample, but feel that it is of
sufficient interest to be worth describing.  A byproduct of this work 
is a slight generalization of Albert's characterization of cyclic algebras 
of prime degree.  See Theorem~\ref{prop.3.6} below.

We will be working here with cyclic algebras and symbol algebras. 
Our notation for these is as follows:  If $T$ is a cyclic Galois field 
extension of a field $K$ of degree $m$ with $\mathcal G(T/K) = 
\langle \tau\rangle$ and $a\in K^*$, we write $(T/K, \tau, a)$ for the 
$m^2$-dimensional cyclic $K$-algebra $\bigoplus\limits_{i=0}
^{m-1}Tx^i$, where $x^m = a$ and $xcx^{-1} = \tau(c)$ for $c\in T$ . If
$\mu_m\subseteq K$ (so $\charr(K)\nmid m$) and $\zeta\in 
\mu_m^*$ and $a,b\in K^*$ we write $(a,b;K)_\zeta$ for the
$m^2$-dimensional symbol algebra over $K$ with generators $i,j$ and
relations $i^m = a$, $j^m = b$, and $ij = \zeta ji$.  For any integer $k$ with 
$\gcd (k,m) = 1$, note the isomorphism 
\begin{equation}
(a,b;K)_\zeta\ \cong \ (a^k,b;K)_{\zeta^k}\,.
\label{eq.3.0}
\end{equation}
For, if $i,j$ are standard generators of $(a,b;K)_\zeta$ as above, 
then $i^k,j$ also generate $(a,b;K)_\zeta$, and they satisfy the 
relations on generators of $(a^k,b;K)_{\zeta^k}$.

Throughout this section $F$ is a field with $\charr (F) \not= p$ and
$\mu_p \nsubseteq F$, and $L = F ( \mu_p)$.  Also, $\mathcal G$,
$J$, and $H$ are as defined in the Introduction.  Let $\Br (F)$ denote
the Brauer group of equivalence classes of central simple
$F$-algebras, and for a field $K \supseteq F$, let $\Br (K/F)$ denote
the kernel of the scalar extension map $\Br (F) \to \Br (K)$; let $_{p^n}
\Br (F)$ (resp.\ $_{p^n} \Br (K/F)$) denote the $p^n$-torsion subgroup of
$\Br (F)$ (resp.\ $\Br (K/F)$).  We have the standard isomorphism
$\Br (L(p)/F) \cong H^2 ( \mathcal G , L(p)^*)$.  The long exact
cohomology sequence arising from the short exact sequence of
$\mathcal G$-modules \eqref{eq.1.14a} above,
 together with the cohomological Hilbert~$90$ \cite[p.~124, Prop.~3]{CF}, 
which says that
$H^1 ( \mathcal G, L(p)^*) = 0$, yields the familiar fact that $_{p^n} \Br
(L(p)/F) \cong H^2 ( \mathcal G, \mu_{p^n})$. 
The Merkurjev-Suslin Theorem (see \cite[Th.~11.5]{MS} or
\cite[p.~149, Th.~8.5]{Sr}) says (as $\mu_p \in L(p)$) that ${}_pBr(L(p))$ is 
generated by symbol algebras of degree $p$.  Since $L(p)$ has no cyclic
field extensions of degree $p$, we must have 
 ${}_p \Br (L(p))= 0$, so ${}_{p^n} \Br (L(p))= 0$; this yields
\begin{equation}
_{p^n}\Br (F)
\cong H^2 ( \mathcal G, \mu_{p^n} ) , \qquad
\mbox{and likewise,} \qquad
_{p^n} \Br (L) \cong H^2 (G_L , \mu_{p^n} ) .
\label{eq.3.1}
\end{equation}

\begin{lemma} \label{lem.3.1}
$_{p^n} \Br (F) \cong (_{p^n} \Br (L))^H$ where $H = \mathcal G (L/F)$. 
Moreover, this isomorphism preserves the Schur index.
\end{lemma}

\begin{proof}
Because $G_L$ is a (in fact, the unique) $p$-Sylow subgroup of
$\mathcal G$ and $\mu_{p^n}$ is $p$-primary torsion, the restriction map
$\res
\!\! : H^2 ( \mathcal G , \mu_{p^n} ) \to H^2 (G_L , \mu_{p^n} )^H$ is an
isomorphism.  (This can be seen
by considering $\cor \circ \res$ and
$\res \circ \cor$, as noted in the proof of Lemma~\ref{lem.2.2}.) 
The isomorphism in the lemma now follows by \eqref{eq.3.1}, since
the second isomorphism in \eqref{eq.3.1} is a $\mathcal G$-module
isomorphism.  (See e.g.~\cite[p.~50]{D} for the action of $\mathcal G$
on $\Br (L)$.)

The restriction map in cohomology corresponds to the scalar
extension  map $_{p^n} \Br (F) \to (_{p^n} \Br (L))^H$.  For a central simple
$F$-algebra $A$, if $A \cong M_t (D)$, i.e., $t \times t$ matrices over a
division ring $D$, then by definition the Schur index of $A$ is $\ind
(A) = \sqrt{\dim_F (D)}$.  If $[A] \in {}_{p^n} \Br (F)$, then  $\ind (A)$ is a
power of $p$, so $\dim_F (D)$ is a $p$-power.  Then, $D \otimes_F L$
is a division algebra, since $\gcd\big  ([ L : F ], \dim_F (D) \big) = 1$ 
(cf.~\cite[p.~67, Cor.~8]{D}).  
Thus, $\ind (A \otimes_F L) = \ind (D \otimes_F L)
= \ind (D) = \ind (A)$, as desired. \hfill $\square$
\end{proof}

\begin{corollary}\label{cor.3.2}
For $J = F(p) ( \mu_p)$, if $(_p \Br (J))^{\mathcal G (J / F(p))} \not=
0$, then there exist division algebras of degree $p$ over $F(p)$ which
are not cyclic algebras.
\end{corollary}

\begin{proof}
By Lemma~\ref{lem.3.1}, with $F(p)$ replacing $F$, if $(_p \Br
(J))^{\mathcal G (J/F(p))} \not= 0$, then $_p \Br (F (p)) \not= 0$.  By
a theorem of Merkurjev \cite[Th.~2]{M83}, the group $_p \Br (F (p))
$ is generated by algebras of degree $p$.  No such algebra can be a
cyclic algebra, since $F(p)$ has no cyclic field extensions of
degree $p$ (see Cor.~\ref{cor.1.2} above).  \hfill $\square$
\end{proof}

We will give examples in \S4 of fields $F$ such that $_p \Br (J)
\not= 0$, but the far more difficult question of nontriviality of $(
_p \Br (J))^{\mathcal G (J/F(p))}$ remains unsettled.

Lemma~\ref{lem.3.1} suggests a further possibility:  There may be
a central simple division algebra $A$ over $L$ of degree $p$ with
$[A] \in \Br (L)^H$, such that $A$ is a cyclic algebra over $L$, but
the inverse image of $A$ in $_p \Br (F)$ is not a cyclic algebra.  This
possibility becomes more plausible when we recall that the cyclic
field extensions of $F$ of degree $p$ correspond only to
a portion of those of
$L$, cf.\ Prop.~\ref{prop.1.4}.

We can put this in sharper focus using the $H$-eigendecomposition
of $_p \Br (L)$ and $L^*/L^{*p}$, where $H = \mathcal G(L/F)$.
 Let $\chi_1, \ldots , \chi_s$ be the
distinct characters $\chi_i \!\! : H \to (\mathbb Z /p \mathbb Z)^*$,
with $\chi_1$ the trivial character;
let $\alpha\!\!:H \to  (\mathbb Z /p \mathbb Z)^*$ 
be the cyclotomic character, as in \eqref{eq.1.9}
above.
  Since $H$ acts on the $p$-torsion
abelian group $_p \Br (L)$, we have $_p \Br (L) =
\bigoplus\limits^s_{i=1} {}_p \Br (L)^{( \chi_i )}$ as in \eqref{eq.1.4}
above.  Lemma~\ref{lem.3.1} shows that $_p \Br (F) \cong {}_p \Br
(L)^{( \chi_1 ) }$.  A central simple $L$-algebra $A$ of degree $p$
which is a cyclic algebra is a symbol algebra, $A \cong (a, b;
L)_\omega$, where $\omega\in \mu_p^*$.  For 
$\sigma\in H$, we have 
 $\sigma [(a,b; L)_\omega] = [(\sigma(a), \sigma (b); L)_{\sigma(
 \omega )}]$.  Note the complication introduced because 
$\sigma$ acts on $\omega$, as well as on $a$ and $b$.
If $[a]$ and $[b]$ lie in eigencomponents of $L^* /
L^{* p}$, then $[A]$ lies in an eigencomponent of $_p \Br (L)$, as we
now  describe.  
The next lemma appears in \cite{M83}.  We include a short proof 
for the convenience of the reader.

\begin{lemma}\label{prop.3.3}
%
%
Let $\chi, \psi$ be characters:  $H \to ( \mathbb Z / p
 \mathbb Z)^*$.   If $[a] \in (L^* /
L^{*p})^{(
\chi )}$ and $[b] \in (L^* / L^{*p})^{(\psi)}$, then $[(a, b; L)_\omega ] \in
{}_p \Br (L)^{(\chi \psi \alpha^{-1})}$.
\end{lemma}

\begin{proof}  
%
%
Recall \cite[p.~80, Lemma 4]{D} that the symbol algebra 
$(a,b; L)_\omega$ depends up to isomorphism
only on the classes $[a]$ and $[b]$ of $a$ and $b$ in $L^*/L^{*p}$. With 
(\ref{eq.3.0}) above, this yields that for any $\sigma\in H$,
$$
\sigma(a,b; L)_\omega\ \cong \ (\sigma(a),\sigma(b);L)_{\sigma(
\omega)} \ \cong  \ (a^{\chi(\sigma)},b^{\psi(\sigma)};L)_{\omega^
{\alpha(\sigma)}}\ \cong \ 
(a^{\chi(\sigma) \alpha(\sigma)^{-1}}, b^{\psi(\sigma});L)
_\omega\, .
$$
Because the class $[(a,b; L)_\omega]$ in ${}_p\Br(L)$ is bimultiplicative
in $a$ and $b$ \cite[p.~80, Lemma 4]{D}, this yields
$\sigma[(a,b; L)_\omega] =[(a,b; L)_\omega]^{\chi(\sigma)\psi(\sigma)
\alpha(\sigma)^{-1}}$  in ${}_p\Br(L)$, as desired.  
 \hfill $\square$
\end{proof}

\begin{proposition}\label{prop.3.4}
Let $\chi \! : H \to ( \mathbb Z/p \mathbb Z)^*$ be a character. 
Take any $a, b \in L^*$ with $[a] \in (L^* / L^{*p} )^{(\chi )}$ and $[b]
\in (L^* / L^{*p})^{(  \alpha
 \chi^{-1})}$, and let $A = (a, b;
L)_\omega$.  Then, there is  central simple algebra $B$ of degree
~$p$  over $F$ with $B \otimes_F L \cong A$.  Furthermore, $B$ is a
cyclic algebra iff there exist $a^\prime , b^\prime \in L^*$ with
$[a^\prime]
\in (L^* / L^{*p} )^{(  \alpha)}$ and $[b^\prime ] \in (L^* /
L^{*p})^H$ such that $A \cong (a^\prime , b^\prime; L)_\omega$.
\end{proposition}

\begin{proof}
By Lemma~\ref{prop.3.3} and Lemma~\ref{lem.3.1}, $[A] \in (_p \Br
(L))^{( \alpha^{-1} \chi ( \alpha \chi^{-1}))} = (_p
\Br (L))^H = \im (_p \Br (F))$.  Because the scalar extension map $_p
\Br (F) \to {}_p \Br (L)$ is index-preserving (see
Lemma~\ref{lem.3.1}), there is a central simple $F$-algebra $B$ 
of degree $p$ with
$B \otimes_F L \cong A$.  Suppose $B$ is a cyclic algebra, say $B
\cong (S/F, \sigma , c)$ where $S$ is a cyclic field extension of $F$
of degree $p$,
 $\mathcal G (S/F) = \langle \sigma\rangle$, and $c \in F^*$.  
 Let $T = S \cdot L$ which
is a cyclic field extension of $L$ of degree $p$, and let
$\sigma^\prime \in \mathcal G (T/L)$ be the generator such that
$\sigma^\prime |_S = \sigma$.  We have $T = L ( \sqrt[p] {a^\prime}
)$ for some $a^\prime \in L^*$, and $a^\prime$ can be chosen so that
$\sigma^\prime ( \sqrt[p]{a^\prime} ) = \omega^{-1}
\sqrt[p]{a^\prime}$.  By Prop.~\ref{prop.1.7}, 
$[a^\prime] \in (L^* / L^{*p})^{(
 \alpha)}$, while $[c] \in F^* / F^{*p} \cong (L^* /
L^{*p})^H$.  Thus, we have $A \cong B \otimes_F L \cong (T/L,
\sigma^\prime , c ) \cong (a^\prime, c; L)_\omega$, as desired. 

Conversely, suppose $A\cong (a',b';L)_\omega$, as in the Prop.
Since ${a'}\in (L^*/L^{*p})^{(\alpha)}$, we know by 
Prop.~\ref{prop.1.7} that there is a cyclic field extension $S$ of $F$ of
degree
$p$, such that $S \cdot L = L(\sqrt[p]{a'})$.  Let $\sigma'$ be the 
generator of $\mathcal G(S\cdot L/L)$ such that 
$\sigma'(\sqrt[p]{a'}) = \omega ^{-1}\sqrt[p]{a'}$, and let 
$\sigma = \sigma'|_S$.  Since $[b'] \in (L^*/L^{*p})^H
\cong F^*/F^{*p}$, there is $c\in F^*$ with $[c] = [b']$ in 
$L^*/L^{*p}$.  Then, $B\otimes_F L \cong A\cong (S/F,
\sigma ,c)\otimes_F L$, so $B\cong (S/F,\sigma, c)$
since the map $\Br(F) \to \Br(L)$ is injective by 
Lemma~\ref{lem.3.1}.  \hfill $\square$
\end{proof}

\begin{remark}\label{rem.3.5}
Prop.~\ref{prop.3.4} suggests a way of obtaining a noncyclic algebra
of degree $p$ over $F$, but we must necessarily start with a 
character $\chi\! :
H \to (\mathbb Z/p \mathbb Z)^*$ different from $ \alpha$ and
the  trivial character $\chi_1$.  We would need 
$[a]\in (L^*/L^{*p})^{(\chi)}$ and $[b]\in (L^*/L^{*p})^{ \alpha \chi^{-1})}$
such that $A = (a,b;L)_\omega$ is a division algebra, but $A$ is not 
expressible as $(a',b';L)_\omega$ for any
$[a^\prime]
\in (L^* / L^{*p} )^{( \alpha)}$ and $[b^\prime ] \in (L^* /
L^{*p})^H$. If $[L:F]\le 2$ then there are not enough different characters, 
and the Prop.~is of no help.  In this connection, recall Merkurjev's result
\cite[Th.~1, Lemma 2]{M83}
that if $[L:F]\le 3$, then ${}_p\Br(F)$
is generated by cyclic algebras of degree $p$.
\end{remark}

The approach in Prop.~\ref{prop.3.4} leads to a slight generalization 
of Albert's characterization of cyclic algebras of prime degree.  This
theorem has recently been  proved independently by 
U.~Vishne, see \cite[Th.~11.4]{V}.

\begin{theorem}\label{prop.3.6}
Suppose $p\nmid [F(\mu_{p^n}):F]$.  Let $D$ be a division algebra
of degree $p^n$ over $F$.  Then, $D$ is a cyclic algebra
over $F$ if and only if there is a $\gamma\in D$ with 
$\gamma^{p^n}\in F^* - F^{*p}$.
\end{theorem}

\begin{proof}
Suppose first that $D$ is a cyclic algebra, say
$D\cong (C/F, \sigma, b)$.  Then, there is $\gamma\in D$
with $\gamma^{p^n} = b$.  If $b\in F^{*p}$, say $b = d^p$,
then for $\delta = \gamma^{p^{n-1}}d^{-1}$ we have
$\delta\in D - F$ and $\delta ^p = 1$.  So, 
$1 < [F(\delta):F] <p$, contradicting 
$[F(\delta):F]\, \big | \,  \dim_F(D)$.  Hence, $b\in F^*-F^{*p}$.

Conversely, suppose there is $\gamma\in D$ with 
$\gamma^{p^n}\in  F^*-F^{*p}$, say $\gamma^{p^n} = c$.
Let $M = F(\mu_{p^n})$.  The assumption that 
$p\nmid [M:F]$ implies that $M= F(\mu_p)$ (see 
Lemma~\ref{lem.1.10}).  
Let $E = D\otimes_F M$. Since $E$  contains the cyclic
Galois field extension $M(\gamma)$ of degree $p^n$ over 
$M$, this $E$ must be a cyclic $M$-algebra; hence, 
$E\cong (a,c;M)_\omega$ for some $a\in M^*$ and 
$\omega \in \mu_{p^n}^*$.  Let 
$\chi_1, \ldots, \chi_s$ be the distinct characters mapping 
$H= \mathcal G(M/F)\rightarrow (\mathbb Z/p^n\mathbb Z)^*$, 
with $\chi_1$ the trivial
character, and let $ \alpha$ be the cyclotomic character (see 
(\ref{eq.1.9})). So
we have, as
$p\nmid [M:F]$,  the eigendecompositions
\begin{equation}
M^*\big /M^{*p^n}\cong \textstyle\prod \limits _{i=1}^s
(M^*\big/M^{*p^n})^{(\chi_i)} \quad\text{ and } \quad
{}_{p^n}\Br(M) \cong \textstyle
\bigoplus\limits_{i=1}^s({}_{p^n}\Br(M))^{(\chi_i)}
\label{eq.3.5}
\end{equation}
by (\ref{eq.1.4}) above.  Write $[a] = \prod\limits_{i = 1}^s[a_i]$ in 
$M^* \big / M^{*p^n}$, where $[a_i]\in (M^*\big /M^{*p^n})^{(\chi_i)}$.
Then, in ${}_{p^n}\Br(M)$, we have $E\cong (a,c;M)_\omega
\sim \bigotimes\limits_{i=1}^s(a_i,c;M)_\omega$.
Also,  $[c]\in (M^*\big /M^{*p^n})^{(\chi_1)}$ as $c\in F^*$;
so each
$(a_i, c;M)_\omega\in {}_{p^n}Br(M)^{(\chi_i \alpha^{-1})}$
by the $p^n$ analogue to Lemma~\ref{prop.3.3}.
Thus, each 
$(a_i,c;M)_\omega$ lies in a different direct summand of
${}_{p^n}\Br(M)$ in the eigendecomposition of 
\eqref{eq.3.5}.  Since 
$[E] \in {}_{p^n}\Br(M)^H ={}_{p^n}\Br(M)^{(\chi_1)}$, we must have 
$E\sim (a_j,c;M)_\omega$ in ${}_{p^n}\Br(M)$, where 
$\chi_j \alpha^{-1} = \chi_1$, i.e., $\chi_j = \alpha$;
dimension count shows that $E\cong (a_j,c;M)_\omega$.
But, since $[a_j]\in (M^*\big /M^{*p^n})^{(\alpha)}$,
the field extension $M(\root p^n\of {a_j}) = S\cdot M$
for some cyclic field extension $S$ of $F$ of degree $p^n$,
by Prop.~\ref{prop.1.7}.  Then, $E\cong(S/F, \tau, c)\otimes_F M$ for 
some generator $\tau$ of $\mathcal G(S/F)$.  Since the map 
${}_{p^n}\Br(F) \to {}_{p^n}\Br(M)$ is injective by Lemma~\ref{lem.3.1}
(as $M = L$),  we have 
$D\cong (S/F, \tau, c)$, as desired.
\hfill $\square$
\end{proof}

\begin{remark}\label{rem.3.7}
Albert's result is the $n = 1$ case of Theorem~\ref{prop.3.6} (see
\cite[Th.~5]{A34} or \cite[p.~177, Th.~4]{A$_2$}, for which
the  condition $p\nmid [F(\mu_p):F]$ always holds.  Our proof 
of Theorem~\ref{prop.3.6} is similar to Albert's, though Albert 
used different terminology, which somewhat veiled his use of 
eigendecompositions.  The theorem is false without the 
assumption that $p\nmid [F(\mu_{p^n}):F]$.  Albert gave in 
\cite{A38} a counterexample with $p^n = 4$, and there are presumably
examples  with odd $p$ also.
\end{remark}

\vskip 0.4truein

\setcounter{equation}{0}
\setcounter{theorem}{0}
\setcounter{section}{4}

\section*{\large 4\ \ Valuations on $J$}

\ \indent
As usual, let $L = F(\mu_p)$, with $L\ne F$ (so $p\ne 2$), 
let $H = \mathcal G(L/F)$, and let 
$J = F(p)(\mu_p)$.  In this section we will look at some of the 
mod $p$ arithmetic of $J$ in order to investigate ${}_p\Br(J)$.
This is motivated by the question discussed in 
\S 3, whether ${}_p\Br(J)^H$ can be nontrivial.
We will use valuation theory, which
is  sometimes a useful tool in verifying that central simple algebras are
division algebras.

\begin{remark}\label{remark4.1} 
Take any field $K$ with $L\subseteq K \subseteq J$ and 
$[K:L] = p$.  Then,  by 
 Albert's theorem (see Prop.~\ref{prop.1.7}
above), $K = L(\root p\of c)$, for $c\in L^*$ with    
$[c]\in \eig L \alpha$, where 
$\alpha\!\!:H\to (\zz/p\zz)^*$
is the cyclotomic
character, as in  \eqref{eq.1.9} above (with $n = 1$).
 By Kummer theory, the map 
$L^*/L^{*p} \to K^*/K^{*p}$ has kernel 
$\langle [c]\rangle\subseteq \eig L \alpha$.  Consequently, for any 
other character $\chi\!\!:H\to (\zz/p\zz)^*$, $\chi\ne \alpha$, the map 
$\eig L \chi \to \eig K \chi$ is injective.  It follows by iteration and passage
to the direct limit that the map $\eig L \chi \to \eig J \chi$ is injective
for each $\chi \ne \alpha$.  Thus, $\eig J \chi$ can be nonzero for each 
$\chi \ne \alpha$, though necessarily $\eig J \alpha = 1$ by 
Albert's theorem, 
as $F(p)$ has no Galois extensions of degree $p$.  It is a more difficult
question when or whether division algebras of degree $p$ over 
$L$ remain division algebras after scalar extension to $J$.  We will use
valuation theory to make some inroads into this question.
\end{remark}  

We will use the following notation: 
Suppose $K$ is a field and $W$
is a valuation ring of $K$.  (This means, in particular, 
that $W$ has quotient field $K$). Let 
$M_W$ denote the maximal ideal of $W$;  let $\ov W = W/M_W$, the 
residue field of $W$;  and let $\vg W$ denote the value group of 
$W$ (written additively).  For a field $K'\supseteq K$, an
extension of $W$ to $K'$ is a valuation ring $W'$ of $K'$ such that 
$W'\cap K = W$.  

\begin{example}\label{laurent}
Let $k$ be any field with $\charr(k) \ne p$ and 
$\mu_p\not\subseteq k$. Let $F$ be the twice iterated
Laurent power series field $F = k((x))((y))$.  Then, 
$F$ has the Henselian valuation ring 
$V = k[[x]]  +yk((x))[[y]]$, with $\ov V \cong k$.  If 
$v\!\!: F^* \to \vg V$ is the associated valuation, then 
$\vg V = \zz\times\zz$, with right-to-left lexicographical 
ordering,
 with $v(x) = (1,0)$ and 
$v(y) = (0,1)$.  
Then, $F(p) = k(p)((x))((y))$ and 
$J = k(p)(\mu_p)((x))((y))$, 
while $L(p) =\bigcup\limits_{i = 1}^\infty k(\mu_p)(p)((x^{1/{p^i}}))
((y^{1/{p^i}}))$.  
(The description of $F(p)$ and $J$ follows 
from Th.~\ref{thfptoj}  below, but can be seen more directly using the
fact that since $\mu_p\not \subseteq k = \ov V$ there is no field
extension of $F$ of degree a power of $p$ which is totally ramified 
with respect to $V$, cf.~\cite[pp.~161--162, (20.11)]{E}
or, more explicitly, \cite[Cor.~2.4]{JWncp} .)  The unique extension of 
$V$ to $J$ is $Z = k(p)(\mu_p)[[x]] + yk(p)(\mu_p)((x))[[y]]$
with $\vg Z = \vg V = \zz\times \zz$; let $z\!\!:J^* \to \vg Z$
be the associated valuation.  For any $\omega\in \mu_p^*$, 
let $D =(x, y;J)_\omega$ (see \S 3 for the notation).  Because 
the images of $z(x)$ and $z(y)$ in $\vg Z/p\vg Z$ are 
$\zz/p\zz$-independent, we know by \cite[Cor.~2.6]{JWncp} that $z$
extends to  a valuation on $D$; so, in particular, $D$ is a division ring. 
Thus, ${}_p\Br(J)$ is nontrivial. 
Since $x$ and $y$ are $H$-stable, it is tempting to think  that $[D]$
should be $H$-stable.  But, in fact, $[D]$ lies in 
${}_p\Br(J)^{(\alpha^{-1})}$ but not in ${}_p\Br(J)^H$ because of 
the nontrivial action of $H$ on~$\mu_p$. 
\end{example}

We will consider valuation rings on $J$ as extensions of ones on $V$.  
For this, let us now fix a valuation ring $V$ of $F$ with $\charr(\ov V)\ne
p$.  Let $W_1, \ldots, W_{\ell}$ be the extensions of $V$ to $L$.  
Let $T = W_1 \cap \ldots \cap W_{\ell}$, which is the integral closure of 
$V$ in $L$.  This notation will be fixed for the rest of this section.
Recall \cite[pp.~95--96, Th.~(13.4)]{E} that the maximal ideals of $T$ are
$N_1,
\ldots N_{\ell}$, where $N_i = M_{W_i} \cap T$,  and that each $W_i$ is the
localization
$T_{N_i}$.

\begin{proposition}\label{propftol}
Let $V$ be a valuation ring of $F$ with $\charr(\ov V) \ne p$.
Let $W_1, W_2, \ldots, W_{\ell}$ be the valuation rings of $L$ 
 extending $V$.  Then, each
 $\ov{W_i} \cong \ov V(\mu_p)$ and
$\vg{W_i} = \vg V$.
Also,  $\ell\,[\ov V(\mu_p):\ov V] = [L:F]$.
\end{proposition}

\begin{proof}
Let $\omega\in \mu_p^*\subseteq L$, and let $f\in F[x]$ be the monic
minimal polynomial of $\omega$ over $F$.  Then, $f\in V[x]$ as 
$\omega$ is integral over $V$, which is integrally closed; 
also, $f\div \sum\limits_{i=0}^{p-1} x^i$ in $F[x]$, and hence
in $V[x]$ by the Division Algorithm, as $f$ is monic.  So, the image
$\ov f$ of $f$ in $\ov V[x]$ divides $\sum\limits_{i=0}^{p-1} x^i$
in $\ov V[x]$.  This shows that the roots of $\ov f$ are all 
primitive $p$-th roots of unity, and $\ov f$ has no repeated roots.
So, if $\ov f = \prod\limits _{i=1}^k g_i$ is the irreducible 
monic factorization  of $\ov f$ in $\ov V[x]$ then the $g_i$ are distinct 
and $\deg(g_i) = [\ov  V(\mu_p): \ov V]$.  Since $fF[x]\cap V[x]  = 
fV[x]$ by the Division Algorithm, we have $V[\omega]\cong
V[x]/fV[x]$, so  
$$
V[\omega]/M_VV[\omega] \  \cong  \ V[x]/(M_V, f)  \ \cong  \ 
\ov V[x]/(\ov f) \ 
\cong  \ \textstyle\bigoplus \limits_{i=1}^k \ov V[x]/(g_i) \, ,
$$
a direct sum of fields.
The inverse images in $V[\omega]$ of the $k$ maximal ideals of 
$V[\omega]/M_VV[\omega]$ are maximal ideals 
$P_1, \ldots, P_k$ of $V[\omega]$
such that each $P_i\cap V = M_V$ and 
$V[\omega]/P_i\cong \ov V[x]/(g_i) \cong \ov V(\mu_p)$.
Because $T$ is integral over $V[\omega]$, for each $P_i$ there is a
maximal ideal 
$N_i$ of $T$ with $N_i\cap V[\omega] = P_i$.  
Then, for $W_i = T_{N_i}$, we 
have $\ov {W_i}\cong T/N_i
\supseteq V[\omega]/P_i \cong \ov V(\mu_p)$.  By the 
Fundamental Inequality, \cite[p.~128, Cor.~(17.8)]{E} or 
\cite[Ch.~VI, \S 8.3, Th.~1]{B}, we have
\begin{eqnarray}\label{fi}
[L:F]&\ge& \textstyle \sum \limits_{i=1}^k[\ov{W_i}:\ov  V]\
 |\vg {W_i} :\vg V| \  \ge \   \sum \limits_{i=1}^k[\ov{W_i}:\ov  V]
\,   \nonumber\\
&\ge&  \textstyle \sum \limits_{i=1}^k[\ov V(\mu_p):\ov V] \  =  \ 
  \sum \limits_{i=1}^k\deg(g_i)  \ = \  \deg(f)  \ = \  [L:F]\ .
\end{eqnarray}
Hence, equality must hold throughout
\eqref{fi}.  
Therefore, each 
$\ov {W_i}= \ov V(\mu_p)$ and $\vg {W_i} = \vg V$, and $k=
[L:F]\big /[\ov V(\mu_p):\ov V]$.  Furthermore,  \eqref{fi} and the 
Fundamental Inequality show that 
$W_1, \ldots, W_k$ are all the extensions of $V$ to $L$; so $\ell
= k$. 
\hfill $\square$ 
\end{proof}

\begin{remark}\label{samel}
Let $S$ be any Galois extension field of $F$ of degree $p$, 
and let $U$ be any extension of $V$ to $S$.  Then, 
$[\ov U:\ov V] \div [S:F] = p$, as $S/F$ is Galois.  Consequently, $\ov U$
and $\ov V(\mu_p)$ are linearly disjoint over $\ov V$, 
and hence $[\ov U(\mu_p):\ov U] = [\ov V(\mu_p):\ov V]$.  It follows
by Prop.~\ref{propftol} applied to $U$ in $S$ in place
of $V$ in $F$ that the number of extensions of $U$ to $S(\mu_p)$
is $\ell$.  Since any field $S'$ with $F\subseteq S'\subseteq F(p)$
and $[S':F]<\infty$ is obtainable from $F$ by a tower of degree $p$ 
Galois extensions (see Prop.~\ref{prop.1.1})
 it follows by iteration that every 
extension of $V$ to $S'$ has exactly $\ell$ extensions to 
$S'(\mu_p)$.  Because this holds for every finite
degree extension $S'$ of $F$ in $F(p)$, it clearly holds for 
every field $S^{\prime\prime}$ with $F\subseteq S^{\prime\prime}
\subseteq F(p)$.
\end{remark}

The   main result of this section describes the residue
field and value group of any  extension of $V$ to $J$.  In
case
$\mu_p
\subseteq \ov V$ (i.e., 
$\ell =  [L:F]$, by Prop.~\ref{propftol}), 
this will require looking at two pieces of $\vg V$.  
For this, let $P$ be the union of all prime ideals $\frak P$ of $V$ such
that 
$V/\frak P$ contains no primitive $p$-th root of unity.  Since the prime
ideals  of $V$ are linearly ordered by inclusion, it is clear that $P$ is a
prime  ideal of $V$ (possibly $P = (0)$), and $P$ is maximal with the
property  that $\mu_p\not\subseteq V/P$.  
(Note also that for every prime ideal $Q\subseteq P$, we have 
$\mu_p\not \subseteq V/Q$.  For, if $\mu_p\subseteq V/Q$,
then $\mu_p\subseteq V/P$, as $\charr(V/P)\ne p$.)
The localization 
$V_P$ of $V$ at $P$ is a valuation ring of $F$ (a \lq\lq coarsening" of
$V$);  let 
$\wi V = V/P$, which is a valuation ring of $\ov{V_P}$.  Recall 
\cite[Ch.~VI, \S 4.3, Remark]{B} that 
there is a canonical short exact sequence of value groups:
\begin{equation}\label{vgses}
0  \ \longrightarrow  \ \vg {\wi V} \  \longrightarrow \  \vg V \ 
\longrightarrow  \ \vg {V_P} \  \longrightarrow \  0
\end{equation}

\begin{theorem}\label{thfptoj}
Let $V$ be a valuation ring of $F$ with 
$\charr(\ov V)\ne p$, and let 
$\ell$ be the number of extensions of
$V$ to $L$.  Let
$Y$be any extension of $V$ to $F(p)$.
Then, $\ov Y \cong \ov V(p)$.
If $\mu_p\not\subseteq \ov V$, then $\vg Y = \vg V$.  If 
$\mu_p \subseteq \ov V$, let $P$ be the prime ideal of $V$ maximal 
such that $\mu_p\not \subseteq V/P$, as above, and let $Q$ be the 
prime ideal of $Y$ with $Q\cap V = P$;  let $\wi Y =  Y/Q$.
Then, $\vg {Y_Q} = \vg {V_P}$, while 
$\vg {\wi Y} = \zz[1/p]\otimes_\zz\vg {\wi V}$.
Furthermore, $Y$ has exactly $\ell$ different 
extensions $Z_1, Z_2, \ldots, Z_{\ell}$ to $J$, 
and each $\ov{Z_i}\cong \ov Y(\mu_p)$ and $\vg{Z_i} = \vg Y$.
\end{theorem} 

Note that in view of the exact sequence like \eqref{vgses} for 
$\vg Y$,  the theorem determines $\vg Y$ completely.  It says that
when we view $\vg Y$ as in the divisible hull $\mathbb Q \otimes
_\zz \vg V$ of $\vg V$, then $\vg Y$ is the subgroup generated by 
$\zz[1/p]\otimes _\zz\vg {\wi V}$ (the $p$-divisible hull of
$\vg {\wi V}$) and $\vg V$.

To prove the theorem we will analyze the range of possibilities
for value groups and residue fields of extensions of $V$ to 
degree $p$ Galois field extensions of $F$.  This will be done
in terms of the corresponding extensions of $L$, where we can
invoke Kummer theory.  To facilitate the analysis, we need
some information on the eigencomponents of induced modules, which is
given in the next proposition.

Let $H = \langle \sigma \rangle$ be a cyclic group of finite order $s$,
and let $\ov H = \langle \si^{m}\rangle$ for some $m\,|\, s$.  Let 
$A$ be any $\ov H$-module, and let $B$ be the induced 
$H$-module, $B = \indh A = \zz[H]\otimes _{\zz[\ov H]} A$. So, as 
abelian groups $B = \bigoplus\limits_{i = 0}^{m - 1}\si^i\otimes A$,
where each $\si^i\otimes A\cong A$.  The left action of $H$ on $B$
arises from the  multiplication action of $H$ on $\zz[H]$.  That is,
\begin {eqnarray}\label{sigmaonind}
\lefteqn{\si\cdot\big( \id\otimes a_0 + \si \otimes a_1 + \ldots+
\si^{m-1}\otimes a_{m-1}\big) \ =} & & \nonumber\\
& & \ \id \otimes \si^{m}
\cdot a_{m-1} + 
\si\otimes a_0 + \si^2\otimes a_1 + \ldots + \si^{m-1}
\otimes a_{m-2} \,.
\end{eqnarray}  

\begin{proposition}\label{eigenind}
With $H = \langle \si\rangle$ and $\ov H = \langle \si^{m}\rangle$
as above, let $A$ be an $\ov
H$-module which is $e$-torsion for some integer $e$.  Let $B = 
\indh A$, as above.  Let $\chi\!\!: H\to (\zz/e\zz)^*$ be any 
character. Then the projection map $\pi\!\!: B\to A$
given by $\sum\limits_{i=0}^ {m-1} \si^i\otimes a_i 
\mapsto a_0$ maps $B^{(\chi)}$ bijectively onto 
$A^{(\chi|_{\ov H})}$, where $\chi|_{\ov H}$ is the restriction
of $\chi$ to $\ov H$.
\end{proposition}

\begin{proof}
Let $b = \sum \limits_{i=0}^ {m-1} \si^i\otimes a_i \in B$.
Note that since $\si^m\in \overline H$,  we have
$\si^m\cdot b = \sum \limits_{i=0}^ {m-1} \si^i\otimes
\si^m( a_i)$. Now, 
 $b \in B^{(\chi)} $ iff $\si\cdot b = \chi(\si) \, b$, iff 
\begin{equation}\label{eigenconditions}
a_0 = \chi(\si)a_1, \ \ a_1 = \chi(\si) a_2, \ \ \ldots, \ \ a_{m-2} = 
\chi(\si) a_{m-1}, \ \text{ and } \ \si^{m}\cdot a_{m-1} = 
\chi(\si)a_0\, .
\end{equation}
If $b\in B^{(\chi)}$, then $\si^m(a_0) = \pi(\si^m\cdot b)
= \pi(\chi(\si)^mb) = \chi(\si^m)(a_0).$
Hence, $a_0\in A^{(\chi|_{\ov H})}$.  Furthermore, if $a_0= 0$
then \eqref{eigenconditions} shows 
that each $a_i = \chi(\si)^{-i}a_0 = 0$; 
so $\pi$ maps $B^{(\chi)}$ injectively 
to $A^{(\chi|_{\ov H})}$.  On the other hand, if we take 
any $a_0 \in A^{(\chi|_{\ov H})}$, then $\si^{m}\cdot a _0
= \chi(\si^{m})a_0 = \chi(\si)^{m}a_0$; so, if we choose 
$a_1 = \chi(\si)^{-1}a_0$, \ldots, $a_i = \chi(\si)^{-i}a_0$,
\ldots, $a_{m - 1} = \chi(\si)^{-(m-1)}a_0$, then 
$\si^{m}\cdot a_{m-1} = \si^{m}\cdot (\chi(\si)^{-(m-1)}a_0) 
= \chi(\si)^{-(m-1)}\si^{m}\cdot a_0 = \chi(\si)a_0$, so the 
equations in \eqref{eigenconditions} are satisfied, 
showing that $a_0\in \pi(B^{(\chi)})$. Thus, 
$\pi\!\! : B^{(\chi)} \to A^{(\chi|_{\ov H})}$ is a bijection.
\hfill $\square$
\end{proof}

We can now prove Theorem \ref{thfptoj}.

\begin{proofof} {\it of Theorem \ref{thfptoj}}.  
It was noted in Remark~\ref{samel} that $Y$ has exactly 
$\ell$ extensions to $J$.  The assertions about $\ov {Z_i}$
and $\vg Z$ follow by applying Prop.~\ref{propftol} to 
$Y$ in $F(p)$ in place of $V$ in $F$.  It remains to 
analyze $\ov Y$ and $\vg Y$.  For this, we look closely
at what can happen with Galois $p$-extensions of $F$.  
These are difficult to get at directly, so we look at the 
corresponding extensions of $L$.

Let us now select and fix one of the $\ell$ extensions of $V$ to $L$; 
call it $W$.  Let $w\!\!: L^* \to \vg W$ be the associated valuation.  
Now, let $c\in L^* - L^{*p}$ with $[c]\in \eig L \alpha$, and let $K = 
L(\root p\of c)$.  Let $S = F(p) \cap K$, which we know by 
Prop.~\ref{prop.1.7}
is a degree $p$ Galois extension of $F$.  (Moreover, all
such Galois extensions of $F$ arise this way.)  Let $R$ be a 
valuation ring of $K$ with $R\cap L = W$;
let $r\!\!: K^*\to \vg R$
be its valuation, and let $U = R\cap S$, which is a valuation ring of 
$S$ with
$U\cap F = V$.  The description of $R$ and $U$ breaks down into three
possible cases:

\medskip

{\it Case I}.  $w(c) \notin p\vg W$.  Then, since 
$r(\root p \of c) = \frac 1p w(c)\in \vg R$, 
the Fundamental Inequality implies that $\vg R = 
\langle \frac 1p w(c)\rangle +\vg W$.  By Prop.~\ref{propftol}
applied to $U$ in $S$  instead of $V$ in $F$, we have 
$\vg U = \vg R =\langle \frac 1p w(c)\rangle +\vg V$. So 
$|\vg U:\vg V| = p = [S:F]$, and the Fundamental Inequality shows that 
$\ov U = \ov V$ and $U$ is the unique extension of $V$ to $S$. 

\medskip 

{\it Case II}.   $w(c) \in p\vg W$.  Then, by modifying $c$ by a
$p$-th power in $L$, we may assume that $w(c)= 0$.  Let 
$\ov c$ be the image of $c$ in $\ov W$.  For this Case II,
assume that $\ov c \notin \ov W^{*p}$.  Then $\ov R$
contains $\ov {\root p \of c} = \root p \of {\ov c}$ which is
not in $\ov W$.  So, the Fundamental Inequality implies that
$\ov R = \ov W(\root p \of {\ov c})$.  Because 
$p = [\ov R: \ov W] \div [\ov R :\ov V]$ but $p\nmid [\ov R:\ov U]$ by 
Prop.~\ref{propftol} applied to $U$ in $S$,  we have $p\div [\ov U:\ov V]$.  
The Fundamental Inequality implies that $[\ov U :\ov V] = p$, 
$\vg U = \vg V$, and $U$ is the unique extension of $V$ to~$S$.  We noted
earlier that $\ov U$ is Galois over $\ov V$.  A comparison of degrees over
$\ov V$ shows that $\ov R = \ov U\cdot \ov W$ so $\ov R$ is abelian 
Galois over $\ov V$.  Thus, $\ov U$ is the unique cyclic Galois extension of 
$\ov V$ of degree $p$ within $\ov R$.  

\medskip

{\it Case III}.  $w(c) \in p\vg W$, so we may assume $w(c) = 0$.
For this Case III, assume that ${\ov c \in \ov W^{*p}}$. We
claim that there are $p$ different valuation rings 
of $K$ extending $W$.  For, consider the subring $W[\root p\of c]$
of $K$.  Since $x^p - c$ is the minimal polynomial of $\root p
\of c$ over $L$, we have $W[\root p \of c] \cong 
{W[x]\big / \big (W[x] \cap (x^p - c) L[x] \big ) }=
W[x] /(x^p - c) W[x]$, where the last equality follows by the Division
Algorithm for monic polynomials in $W[x]$.  Hence, 
$W[\root p \of c]\big / M_W\, W[\root p \of c] 
\cong W[x] \big / \big (M_W, x^p - c\ \big ) 
\cong \ov W[x] /(x^p - \ov c)$.  Because $\ov c\in 
\ov W^{*p}$ and $\mu_p \subseteq \ov W$, $x^p -\ov c$ factors into 
distinct linear terms in $\ov W[x]$, say $x^p - \ov c = 
(x-d_1) \ldots (x-d_p)$.  Then, the Chinese Remainder Theorem
shows that $\ov W[x] /(x^p - \ov c) \cong \bigoplus\limits
_{i=1}^p \ov W[x] / (x-d_i)$.  Because $W[\root p \of c]\big /M_W
W[\root p \of c]$ 
thus has $p$ maximal ideals, $W[\root p \of c]$ has at least $p$
maximal ideals.  Let $C$ be the integral closure of $W$ in $K$.  Since
$C$ is integral over $W[\root p \of c]$, 
 $C$ has at least $p$ different maximal 
ideals, say  $N_1, \ldots, N_p$.
Each localization $R_i = C_{N_i}$ is a
different  valuation ring of $K$ with $R_i \cap L = W$.  The Fundamental
Inequality shows that there must be exactly $p$ of the $R_i$, as claimed.

Now, since $\mathcal G(S/F)$ acts transitively on the valuation rings of 
$S$ extending $V$ \cite[p.105, (14.1)]{E}, 
the number of such extensions is either $1$ or $p$.
There are  at least $p$ extensions of $V$ to $K$ (namely, the $R_i$),
but every extension of $V$ to $S$ has $\ell \le p-1$ extensions to 
$K$ by Prop.~\ref{propftol} applied over $S$.  Hence, there must
be more than one, so exactly $p$ extensions of $V$ to $S$, 
call them $U_1, \ldots, U_p$.  The Fundamental Inequality shows that
each $\ov {U_i} = \ov V$ and $\vg {U_i} = \vg V$.  This
completes Case III.

\medskip
 
We must still see what constraints are imposed by the condition that 
$[c] \in \eig L\alpha$. For this, let $H = \mathcal G(L/F) = 
\langle \si \rangle$,  as usual, and let  
 $\ov H = \{\, \tau\in H\ |\ \tau(W)
= W\,\}$, the decomposition group of $W$ over $V$.  
Because $H$ acts transitively on the set of extensions of $V$ to 
$L$ and there $\ell$ such extensions, $|H:\ov H| = \ell$, 
so $\ov H = \langle \si^{\ell} \rangle$.  
Each $\tau \in \ov H$ maps $W$ to itself, so induces an automorphism
$\ov \tau$ of $\ov W$. Recall 
\cite[p.~147, (19.6)]{E} or \cite[p.~69, Th.~21]{ZS} that the map $\ov H\to
\mathcal G (\ov W/\ov  V)$ given by $\tau \mapsto \ov \tau$ is a group 
epimorphism.  
By Prop.~\ref{propftol} we have $\big |\ov H \big| = 
|H|/\ell = [L:F]/\ell  =  | \mathcal G(\ov W/\ov  V)|$,
and therefore the map $\ov H\to
\mathcal G (\ov W/\ov  V)$ is an isomorphism.
Also, because $\ov \tau$ acts on the $p$-th roots of unity in 
$\ov W$ according to the action of $\tau$ on the $p$-th roots
of unity in $L$, the cyclotomic character
$\ov \alpha$ for $\mathcal G(\ov W/\ov  V)$
corresponds to the restriction $\alpha|_{\ov H}$.

Observe that  the distinct extensions of 
$V$ to $L$ are $\si^i(W)$ for $0\le i\le \ell - 1$.  Each $\vg{\si^i(W)}$ is 
canonically identified with $\vg W$ inside the divisible hull of 
$\vg V$, and for the associated valuation $w_i$ of $\si^i(W)$ we have 
$w_i = w\circ \si^{-i}$.   Likewise, for $0\le i\le \ell -1$ we identify 
$\ov {\si^i(W)}$ with $\ov W$ using the isomorphism $\ov {\si^i}
\!\!:\ov W \to \ov{\si^i(W)}$ induced by $\si^i\!\!: W \to \si^i(W)$.  
So, for $c\in \si^i(W)$, we have $\ov c\in \ov{\si^i(W)}$ corresponds
to $\ov {\si^{-i}(c)}$ in $ \ov W$. 

\medskip

We can now determine $\overline Y$.

View $\ov W^*$ as an $\ov H$-module, where $\tau\in \ov H$ 
acts by $\ov \tau$.  Let $\indh \ov W^*$ be the induced $H$-module 
described before Prop.~\ref{eigenind}, with $m = \ell$.
Recall that $T$ denotes 
the integral closure of $V$ in $L$, so $T = \bigcap\limits_{i=0}^{\ell
-1}\si^i(W)$ \cite[p.~95, Th.~3.3.(b)]{E}.  
Let $\gamma\!\!:T^* \to \indh \ov W^*$ be the map
given by $\gamma(t) = \sum\limits_{i=0}^{\ell - 1}
\si^i\otimes\ov{\si^{-i}(t)}$ (the bar denotes image in $\ov W^*$).
The surjectivity of $\gamma$ is equivalent to the assertion that 
for every $r_0, \ldots, r_{\ell -1}\in \ov W ^*$ there is $t \in T^*$ with 
$\ov{\sigma^{-i}(t)}  = r_i$  in $\ov W$ for each $i$, i.e., $\ov t = \ov
{\sigma^i} (r_i)$ in $\ov{\sigma^i(W)}$.  This holds by the 
Approximation Theorem \cite[p.~79, Th.~(11.14)]{E} or 
\cite[p.~30, Lemma~2]{ZS}.
(For this the valuation rings $\si^0(W) , \ldots, \si^{\ell-1}(W)$
need not be independent, just incomparable. This result 
uses only the Chinese Remainder Theorem applied to $T$.)
Also, since $\si ^{\ell}\cdot \ov {\si^{-(\ell -1)}(t)} 
= \ov{\si^\ell (\si^{-(\ell-1)}(t))} = \ov{\si(t)}$, we have 
$\si\cdot \gamma(t) = \gamma(\si(t))$, so $\gamma$ is an 
$H$-module epimorphism.  Therefore, the corresponding map
$T^*/T^{*p}\to \indh (\ov W^* /\ov W^{*p})$ is an $H$-module
epimorphism.  So, $\eig T \alpha$ maps onto ${\big(\indh  (\ov W^*/
\ov W^{*p})\big)}^{(\alpha)}$,
which by Prop.~\ref{eigenind} projects onto $\eig{\ov W}{\ov
\alpha}$.    That is, for any $a\in \ov W^* - \ov W^{*p}$ such that 
$[a] \in \eig {\ov W}{\ov \alpha}$ there is $t\in T^*$ with 
$[\ov t] = [a]$ in $\ov W^*/\ov W^{*p}$.  If we choose $c = t$, then 
for the resulting $K = L(\root p \of c)$ we are in Case II above, with 
$\ov R = \ov W(\root p \of {\ov t} )= \ov W(\root p \of a)$, and 
$\ov U$ is the degree $p$ Galois extension of $\ov V$ within 
$\ov R$.  Since we can do this for any $[a] \in \eig {\ov W}{\ov \alpha}$,
Prop.~\ref{prop.1.7} shows that every Galois extension of $\ov V$ of 
degree $p$ is realizable as some $\ov U$, and so lies in $\ov Y$.

Now, $F(p)$ is the direct limit of finite towers of Galois extensions
of degree $p$ starting with $F$ (see Prop.~\ref{prop.1.1}).  If $S'$ is the 
top field in such a tower, then $\ov {Y\cap S'}$ is obtained from $\ov  V$
by a succession of Galois extensions of degree $1$ or $p$.  
Hence $\ov {Y\cap S'} \subseteq \ov V(p)$ for each $S'$, and therefore 
$\ov Y \subseteq  \ov V(p)$. But, iteration of the argument in the
preceding paragraph shows that any finite degree extension of 
$\ov V$ within $\ov V(p)$ is obtainable as  $\ov {Y\cap S'}$ for a 
suitably built $S'$.  Hence, $\ov Y = \ov V(p)$, as desired.

\medskip

We now determine $\vg Y$.  

For the trivial $\ov H$-module $\vg W$, we have
the induced $H$-module $\indh \vg W$.  Let $\beta\!\!: L^* \to \indh \vg W$
be the map given by $d\mapsto \sum\limits _{i=0}^ {\ell-1}
\si^i\otimes w(\si^{-i}(d))$.  Since 
$\si^{\ell}\cdot w(\si^{-(\ell-1)}(d)) = w(\si(d))$, as $w\circ \si^\ell = w$
and
$\si^\ell$ acts trivially  on $\vg W$, this $\beta$ is an $H$-module
homomorphism.  By reducing mod $p$ we obtain an 
$H$-module homomorphism $\ov \beta\!\!: L^*/L^{*p} \to 
\indh(\vg W/p\vg W)$. So, for our $c\in L^*$ used to define $K$, 
since $[c] \in \eig L\alpha$, we have $\ov \beta[c] \in (\indh \vg W/p
\vg W)^{(\alpha)}$, so Prop.~\ref{eigenind} shows that 
$w(c) + p\vg W\in (\vg W/p\vg W)^{(\alpha|_{\ov H})}$.  

Suppose first that $\mu_p\not\subseteq \ov V$.  Then, 
$\ell < s = [L:F]$, by Prop.~\ref{propftol}. so $\ov H$, of order 
$s/\ell$, is nontrivial. Since the cyclotomic character $\alpha$ has order 
$s$, its restriction $\alpha|_{\ov H}$  has order $|\ov H|$, so is nontrivial.  
Since $\ov  H$ acts trivially on $\vg W$, it follows that $(\vg W/p\vg
W)^{(\alpha|_{\ov H})} = (0)$.  Now, the only way we could have 
$\vg U$ larger than $\vg V$ is if our $c$ is in Case I above.  But then
we would have $w(c)\notin p\vg W$, yielding a nontrivial 
element in the trivial group $(\vg W/p\vg
W)^{(\alpha|_{\ov H})}$.  Since this cannot occur, we see that Case~I never
arises when $\mu_p \not \subseteq \ov V$.   Therefore, $\vg U = \vg V$
for every degree $p$ Galois extension $S$ of $F$.  It follows by iteration
and passage to the direct limit that $\vg Y = \vg V$, as asserted.

Now suppose instead that $\mu_p\subseteq \ov V$.  Prop.~\ref
{propftol} shows that $\ell = s$, i.e., there are $s$ different extensions 
$W_1, \ldots, W_s$ of $V$ to $L$.  Consider first the extreme case where 
$\mu_p \subseteq V/\frak p$ for each nonzero prime ideal $\frak p$
of $V$.  For any such $\frak p$, the extensions of the localizations
$V_{\frak p}$ to $L$ are the localizations $W_{1 \frak p}, \ldots,
W_{s \frak p}$.  
(Each $W_{i\frak p}$ coincides with the localization of $W_i$ at its prime
ideal lying over $\frak p$.)
Since $\mu_p\subseteq \ov{V_{\frak p}}$,
 which is the quotient field of $V/\frak p$, Prop.~\ref{propftol} applied to 
$V_{\frak p}$ shows that $V_{\frak p}$ has $s$ different extensions to 
$L$.  (The Prop.~applies, as $\charr(\ov {V_{\frak p}}) \ne p$.)
So, $W_{i\frak p} \ne W_{j\frak p}$ for $i\ne j$.  Now, for each $i$,
the rings between $W_i$ and $L$ are the $W_{i\frak p}$ as 
$\frak p$ ranges over the nonzero prime ideals of $V$.  Since
$W_{i\frak p} \ne W_{j\frak  p}$  for $i\ne j$, it follows that the 
valuation rings $W_1, \ldots, W_s$ are pairwise independent, i.e., 
there is no valuation ring of $L$ (smaller than $L$ itself) containing 
both $W_i$ and $W_j$ for any $i\ne j$.  Because of this independence,
the Approximation Theorem (see \cite [p.~80, (11.16)]{E}) applies, and
shows that our map $\beta\!\!: L^* \to \indh \vg W$ is surjective;
so $\ov \beta\!\!: L^*/L^{*p}\to \indh(\vg W/p\vg W)$ is also surjective,
so it is also surjective when restricted to the $\alpha$-eigencomponents.
By Prop.~\ref{eigenind} $(\indh (\vg W/p\vg W))^{(\alpha)}$ projects onto 
$(\vg W/p\vg W)^{(\alpha|_{\ov H})}$, which here is all of $\vg W/p\vg W$
since $|\ov H| = 1$ as $\ell = s$.  This means that for any $\varepsilon
\in \vg W - p\vg W$ there is $c \in L^*$ such that $[c] \in \eig L\alpha$
and $w(c) \equiv \varepsilon \ (\text{\sl {mod}} \ p\vg W)$.  If we let 
$K = L(\root p \of c)$ for this choice of $c$, then
we are in Case I above, which shows that $\vg U = \vg R = 
\langle \frac 1p \varepsilon \rangle + \vg W$.  Since this is 
true for any $\varepsilon \in  \vg W - p\vg W$, it follows by iteration and 
passage to the direct limit that $\vg Y = \varinjlim \frac1{p^n}\vg V
= \zz[1/p] \otimes _\zz\vg V$.  This is what is asserted in the 
theorem, since in the extreme case we are now considering
$P = (0)$, so $\wi V = V$ and $\wi Y = Y$.

We handle the general situation by combining the cases previously 
considered.  Suppose $\mu_p\subseteq \ov V$.  For the 
prime ideal $P$ defined in the theorem, we have $\mu_p \not \subseteq
\ov {V_P}$, which is the quotient field of $V/P$.  
Now, $Y_Q$ is an extension of $V_P$ to $F(p)$.  
Since $\mu_p\not\subseteq \ov {V_P}$ and $\charr(\ov {V_{\frak p}})
\ne p$, by applying to $V_{\frak p}$
the argument given 
previously for $V$ we obtain 
$\vg {Y_Q} = \vg {V_P}$, as desired.  Furthermore, 
$\ov {Y_Q}\cong \ov{V_P}(p)$.  Thus, $\wi Y = Y/Q$ can be viewed as 
an extension of $\wi V = V/P$ from $\ov {V_P}$ to $\ov {V_P}(p)$.
By the choice of $P$, the extreme case considered in the 
previous paragraph applies to $\wi V$.  
Hence, $\vg {\wi Y}  = \zz[1/p]\otimes _\zz\vg {\wi V}$. 
\hfil $\square$

\end{proofof}

\begin{example}\label{mixedex} 
Let $F_0 = \mathbb Q(x,y)$, the rational function field in 
two variables over $\mathbb Q$.  Let $V_0 = 
\mathbb Q[x]_{(x)} + y\mathbb Q(x)[y]_{(y)}$.
Here, we are localizing first with respect to the prime ideal 
$(x) $ of $\mathbb Q[x]$,  and second with respect to the 
prime ideal $(y)$ of $\mathbb Q(x)[y]$.
Then, $V_0$ is a valuation 
ring of $F_0$ with  $\ov{V_0} \cong \mathbb Q$ and $\vg{V_0} = 
\zz\times \zz$.  If $v_0\!\!: F_0^* \to \vg{V_0}$ is the associated
valuation, then $v_0(x) = (1,0)$ and $v_0(y) = (0,1)$.  Note that 
$V_0$ is the intersection with $F_0$ of the standard Henselian 
valuation ring on $\mathbb Q((x))((y))$ described in Ex.~\ref{laurent}
above.  For any odd prime $p$,
let $F = F_0(\root p \of {1+x})$.  
To see how $V_0$ extends to $F$, let $T$ be the integral closure of 
$V_0$ in $F$, and let $S = V_0[\root p \of {1+x}] \subseteq T$.
Since $S\cong V_0[t]\big/ \big (t^p - (1+x)\big )$, we have 
$$
S/M_{V_0}S \cong \ov{V_0}[t]\big/\big (t^p-(1+\ov x)\big )
\cong \mathbb Q[t]/(t^p-1) \cong \mathbb Q[t]/(t-1)\oplus 
\mathbb Q[t]/(t^{p-1} + \ldots +1)
\cong \mathbb Q \oplus \mathbb Q(\mu_p)\,  .
$$
So, $T$, being integral over $S$, has at least two maximal ideals
$N_1$ and $N_2$, with $\mathbb Q\subseteq T/N_1$ and 
$\mathbb Q(\mu_p)\subseteq N_2$.
The Fundamental Inequality shows that for the extensions $V_i = 
T_{N_i}$ of $V_0$ to $F$, we have 
$\ov {V_1} \cong \mathbb Q$, 
$\ov {V_2} \cong \mathbb Q(\mu_p)$, and $\vg {V_1} = 
\vg {V_2} = \vg{V_0} =\zz\times\zz$.  
Furthermore, $V_1$ and $V_2$ are the only extensions of $V_0$ to $F$.
If $Y_i$ is any extension of 
$V_i$ to $F(p)$,  then Th.~\ref{thfptoj} shows that 
$\ov {Y_1} \cong \mathbb Q(p)$ and $\vg{Y_1} = \zz\times\zz$.
Let $\frak p$ be the prime ideal $yV_2$.  Then, 
$V_2/\frak p \cong \mathbb Q(x)(\root p \of {1+x})$, which does not 
contain~$\mu_p$.  So, $\frak p$ is the prime ideal $P$ of 
Th.~\ref{thfptoj} for $V_2$.  Since $\vg{V_2/\frak p} = 
\zz\times 0$, 
Th.~\ref{thfptoj}
shows that 
$\vg{Y_2} = \zz[1/p] \times \zz$
while $\ov{Y_2} \cong \mathbb Q(\mu_p)(p)$.
\end{example}

\begin{remark}
We had hoped to use valuation theory to construct an example
of a nonsplit algebra of degree $p$ in ${}_p\Br(J)^H$.  However, we will
now show why Th.~\ref{thfptoj} does not help in this.  Let $V$ be a valuation
ring of $F$ with $\charr(\ov V)\ne p$, let $W$ be an extension of 
$V$ to $L$ with associated valuation $w\!\!: L^* \to \vg W$, and let
$Z$ be an extension of  $W$ to $J$.  There are three types of symbol
algebras $A = (a,b;L)_\omega$ (with $a,b\in L^*$ and $\omega\in 
\mu_p^*$) for which it is known that $w$ extends to a valuation
on $A$, and hence $A$ is a division algebra:   \ (1) $w(a)$ and $w(b)$
map to $\zz/p\zz$-independent elements of $\vg W/p\vg W$.  Then, 
cf. \cite[Cor.~2.6]{JWncp},  the 
valuation ring of $A$ is tame and totally ramified over $W$, with residue
division algebra $\ov  V$ and value group $\langle \frac 1p w(a), \frac 1p
w(b)\rangle + \vg W$.   \ (2) $w(a) \notin p\vg W$ and $w(b) = 0$, and for 
the image $\ov b$ of $b$ in $\ov W$ we have $\ov b \notin \ov W^{*p}$.
Then, cf.~\cite[Cor.~2.9]{JWncp},
 the valuation ring of $A$ is semiramified over $W$, with residue 
division algebra $\ov W(\root p\of{\ov b})$ and value group 
$\langle \frac 1p w(a)\rangle + \vg W$. \  (3) $w(a) = w(b) = 0$ and 
$(\ov a, \ov b ;\ov W)_{\ov \omega}$ is a division ring.  Then,
the 
valuation ring on $A$ is unramified over $V$, with residue algebra 
$(\ov a, \ov b ;\ov W)_{\ov \omega}$ and value group $\vg W$.
For, if $i$ and $j$ are standard generators of $A = (a,b;L)_\omega$, 
then it is easy to check that the map $u\!\!:A - \{0\} \to \vg W$ given
by $u\big (\sum\limits_{r=0}^{p-1}\, \sum \limits_{s=0}^{p-1}
c_{rs}\,  i^rj^s \big)= \text {\sl min} \{w(c_{rs})\,|\,c_{rs}\ne 0\}$
($c_{rs}\in L$) is a valuation on $A$ with the specified
residue algebra and value group.  (The proof is similar to but easier
than the proof of \cite[Th.~2.5]{JWncp}.)
 Since type (3) reduces the problem of obtaining a division algebra
to the same problem over the residue field, it is not helpful for
constructing  examples, and we will not consider this type further.

Suppose we choose $a,b \in L^*$ so that for some character 
$\chi\!\!:H\to \zz/p\zz^*$, we have $[a]\in \eig L\chi$ and $[b]\in 
\eig L{\alpha \chi^{-1}}$.  Then, $[A]\in {}_p\Br(L)^H$ by 
Lemma~\ref{prop.3.3}, so $[A\otimes_LJ] \in {}_p\Br(J)^H$.  but, we will
see that the valuation conditions that assure $A$ is a division 
algebra break down over $J$.  Suppose first that $\mu_p\notin 
\ov V$; so, in the notation of Th.~\ref{thfptoj} and its proof, 
$\ell < [L:F]$ and $\ov H$ is nontrivial.  Suppose $w(a)\notin p\vg W$; then 
as in the proof of Th.~\ref{thfptoj}, Prop.~\ref{eigenind} implies that 
the image of $w(a)$ is nontrivial in ${\vg W/p\vg W}^{(\chi|_{\ov H})}$;  
this forces $\chi|_{\ov H}$ to be trivial, as $\ov H$ acts trivially on 
$\vg W$.  Hence, $\alpha \chi^{-1}|_{\ov H} = \alpha|_{\ov H}$, which is 
nontrivial and is identified with the cyclotomic character
$\ov \alpha$ for $\mathcal G(\ov W/\ov V)$.  The nontriviality of 
$\alpha\chi^{-1}|_{\ov H}$ forces $w(b)\in p\vg W$, so we may assume 
$w(b) = 0$.  If $b\notin \ov W^{*p}$, then $A$ is a division algebra of
type (2).  But, Prop.~\ref{eigenind} implies that $\ov b$ maps to 
$\eig{\ov W}{\alpha \chi^{-1} |_{\ov H}} = \eig {\ov W}{\ov \alpha}$.  Hence,
on passing to  $J$ we find that $\ov b\in \eig {\ov Z}{\ov \alpha}$,
which is trivial as $\ov Z \cong \ov V(p)(\mu_p)$---see 
Remark~\ref{remark4.1}.
This means that $\ov b \in \ov Z^{*p}$, and we have lost the conditions for 
type (2) for $A\otimes _L J$.  Likewise, if $w(b) \notin p\vg W$, then
we are forced to have $w(a) \in p\vg W$, and when we adjust $a$ so that
$w(a) = 0$, the same argument as just given shows that 
$\ov a \in \ov Z^{*p}$.  Thus, we have not been able to obtain a type (1) 
or a type (2) valued division algebra in ${}_p\Br(J)^H$ when $\mu_p
\notin \ov V$. 

Suppose instead that $\mu_p\subseteq \ov V$.  Since $\ov Z = \ov V(p)
(\mu_p) = \ov V(p)$ and $\mu_p \subseteq  \ov V$, $\ov Z^*/
\ov Z^{*p}$ is trivial.  Therefore, we will not obtain any valued division
algebras of degree $p$ of type (2) or type (3) over $J$.  We are left to 
search for type (1) division algebras.  Thus, we may assume that 
$w(a)$ and $w(b)$ are $\zz/p\zz$-independent in $\vg W/p\vg W$.  
Here $\ov H$ is trivial, but choose the prime ideal $P$ of $V$
as in Th.~\ref{thfptoj}, and let $\frak P$ be the prime ideal of $W$ with 
$\frak P\cap V = P$, and $\wi H$ the (nontrivial) decomposition 
group of $W_{\frak P}$ over $V_P$; let $w_{\frak P}$ be the valuation 
of $W_{\frak P}$.   We have an $H$-module homomorphism 
$\wi \gamma\!\!: L^*/L^{*p} \to 
\text {\sl{ind}}_{\wi H\to H}( \vg{W_{\frak P}}\big /p\vg {W_{\frak P}})$
so since $[a]\in \eig L\chi$ we find that $\wi \gamma [a]  \in
{\big ({\text {\sl{ind}}}_{\wi H\to H}( \vg{W_{\frak P}}\big /p\vg {W_{\frak
P}})\big)}^{(\chi)}$.  By Prop.~\ref{eigenind} it follows that 
 $w_{\frak P}(a) \in 
{\big(\vg {W _{\frak P}}\big /p \vg {W_{\frak P}} \big)}^
{(\chi|_{\wi H})}$.  Since $\wi H$ acts trivially on $\vg {W_{\frak P}}$,\
this implies that $w_{\frak P}(a) \in p\vg {W_{\frak P}}$ or 
$\chi|_{\wi H}$ is trivial.   If $w_{\frak P}(a) \in p \vg {W_{\frak P}}$, 
we can modify $a$ by a $p$-th power in $L^*$ to assume that $w_{\frak P}
(a) = 0$; but then, for $\wi W = W/\frak P$ the exact sequence like 
\eqref{vgses} for $\vg W$ shows that $w(a) \in \vg {\wi W}$.  But then,
Th.~\ref{thfptoj} shows that $w(a) \in p \vg Z$, so that $(a,b;J)_{\omega}$
is not a type (1)  valued division algebra over $J$.  On the other hand, if 
$\chi|_{\wi H}$ is trivial, then $\alpha\chi^{-1}|_{\wi H} = \alpha|_{\wi H}$,
which is nontrivial.  Hence, the argument just given for $a$ now shows that 
$w(b)\in p\vg Z$, so again we do not obtain a type (1) valued division
algebra over $J$. 
\end{remark}


\vskip0.4truein

\baselineskip=12pt
\parindent=0pt



Department of Mathematics \\
Middlesex College\\
 University of Western Ontario \\
London, Ontario N6A 5B7\\
Canada

{\em e-mail}:  {\tt minac@uwo.ca}

\bigskip\bigskip

Department of Mathematics, 0112\\
 University of California, San Diego \\
 9500 Gilman Drive \\
 La Jolla, CA 92093-0112 \\
 USA

{\it e-mail}: {\tt arwadsworth@ucsd.edu}

\end{document}